\documentclass[preprint,review,12pt]{elsarticle}
\usepackage{amssymb}
\usepackage{amsfonts}
\usepackage{amsmath}
\usepackage{amsthm}
\usepackage{color}
\usepackage{graphicx}
\usepackage{epsfig,mathrsfs}
\usepackage[bf,SL,BF]{subfigure}
\usepackage{fancyhdr}
\usepackage{CJK}
\usepackage{caption}
\usepackage{wrapfig}
\usepackage{cases}

\newtheorem*{lemma*}{Lemma A}

\numberwithin{equation}{section}

\renewcommand{\eqref}[1]{(\ref{#1})}

 \allowdisplaybreaks

\begin{document}
 \pagenumbering{arabic}

\begin{frontmatter}
\title{{\bf \noindent Second order finite difference approximations for the two-dimensional  time-space Caputo-Riesz fractional  diffusion equation}}
\author{Minghua Chen, Weihua Deng$^{*}$, Yujiang Wu}
\cortext[cor2]{Corresponding author. E-mail: dengwh@lzu.edu.cn.}
\address{School of Mathematics and Statistics,
Lanzhou University, Lanzhou 730000, P. R. China}

\date{}

\begin{abstract}
In this paper, we discuss the time-space Caputo-Riesz fractional  diffusion equation
with variable coefficients  on a finite domain. The finite difference schemes for this equation are provided. We theoretically prove and numerically verify that
the implicit finite difference scheme is unconditionally stable
(the explicit scheme is conditionally stable with the stability condition $\frac{\tau^{\gamma}}{(\Delta x)^{\alpha}} +\frac{\tau^{\gamma}}{(\Delta y)^{\beta}} <C$)
 and 2nd order convergent in space direction,
and $(2-\gamma)$-th order convergent in time direction, where $\gamma \in(0,1]$.

\end{abstract}
\end{frontmatter}

\section{Introduction}
Nowadays, fractional calculus has become popular in both the science and engineering societies. There are several nonequivalent definitions of fractional derivatives \cite{Miller:93,Podlubny:99,Samko:93}. The Caputo derivative is most often used for time fractional derivative, and the Riemann-Liouville derivative and Gr\"{u}nwald-Letnikov derivative, the two fractional derivatives being equivalent if the functions performed are regular enough, are most frequently used for space fractional derivative. Some important progress has been made for numerically solving this kind of fractional PDEs by finite difference methods, e.g., see \cite{Chen:11,Liu:07,Meerschaert:06,Sousa:12,Sun:06,Yuste:06,Zhuang:07}.

Another space fractional derivative having a vast majority of applications is the symmetric fractional derivative, namely the Riesz fractional derivative, e.g., see \cite{Metzler:00,Zaslavsky:02}.  Jiang et al analytically discuss the time-space fractional advection-diffusion equation with Riesz fractional derivative as space fractional derivative \cite{Jiang:12}. Based on the shifted Gr\"{u}nwald approximation strategy and the method of lines,  Yang, Liu, and Turner numerically study the Riesz space fractional PDEs with two different fractional orders $1<\alpha \le 2$ and $0 <\beta <1$ \cite{Yang:10}. The explicit finite difference scheme for fractional Fokker-Planck equation with Riesz fractional derivative is discussed in \cite{Deng:12}. With the desire of obtaining 2nd order convergence in the space discretization, here we further discuss the finite difference approximations for two-dimensional time-space Caputo-Riesz fractional diffusion equation with variable coefficients in a finite domain, namely,
\begin{equation} \label{1.1}
\left\{ \begin{array}
 {l@{\quad} l}
 _0^C\!D_t^{\gamma}u(x,y,t)=c(x,y,t)  \displaystyle \frac{\partial ^{\alpha}u(x,y,t)}{\partial |x|^{\alpha}}
 +d(x,y,t)\frac{\partial ^{\beta}u(x,y,t)}{\partial |y|^{\beta}}+f(x,y,t),\\
  u(x,y,0)=u_0(x,y) ~~~~\, {\rm for}~~~ (x,y) \in \Omega,\\
   u(x,y,t)=B(x,y,t) ~~~ {\rm for}~~~ (x,y,t) \in \partial \Omega \times [0,T],
 \end{array}
 \right.
\end{equation}
in the domain $\Omega=(x_L,x_R) \times (y_L,y_R),\, 0< t \leq T$,
with the orders of the Riesz fractional derivative  $1<\alpha,\beta\leq2$
and the order of the Caputo  fractional operator  $0<\gamma\leq1$; the function $
f(x,y,t)$ is a source term; and the variable coefficients
$c(x,y,t)\geq 0, d(x,y,t)\geq0$. The Riesz fractional derivative  for $n \in \mathbb{N}$, $n-1 < \nu < n$, in a finite interval $x_L \leq x \leq x_R$
is defined as  \cite{Samko:93}
\begin{equation}\label{1.2}
\frac{\partial ^{\nu}u(x,y,t)}{\partial |x|^{\nu}} =-\kappa_{\nu}\left( _{x_L}\!D_x^{\nu}+ _{x}\!D_{x_R}^{\nu} \right)u(x,y,t),
\end{equation}
where the coefficient $\kappa_{\nu}=\frac{1}{2cos(\nu \pi/2)}$, and
\begin{equation}\label{1.4}
 _{x_L}D_x^{\nu}u(x,y,t)=
\frac{1}{\Gamma(n-\nu)} \displaystyle \frac{\partial^n}{\partial x^n}
 \int_{x_L}\nolimits^x{\left(x-\xi\right)^{n-\nu-1}}{u(\xi,y,t)}d\xi,
\end{equation}

\begin{equation}\label{1.5}
 _{x}D_{x_R}^{\nu}u(x,y,t)=
 \frac{(-1)^n}{\Gamma(n-\nu)}\frac{\partial^n}{\partial x^n}
\int_{x}\nolimits^{x_R}{\left(\xi-x\right)^{n-\nu-1}}{u(\xi,y,t)}d\xi,
\end{equation}
are the left and right Riemann-Liouville space fractional
derivatives, respectively. The Caputo fractional derivative of order $\gamma \in (0,1]$ is defined by \cite{Miller:93,Podlubny:99}
\begin{equation}\label{1.6}
 _0^{C}\!D_{t}^{\gamma}u(x,y,t)=\left\{ \begin{array}
 {l@{\quad} l}
\displaystyle \frac{1}{\Gamma(1-\gamma)} \int_{0}^t \frac{\partial u(x,y,{\eta})}{\partial {\eta}} {(t-\eta)^{-\gamma}}d\eta,& 0<\gamma <1,\\
 \cr\noalign{\vskip 0 mm}
\displaystyle \frac{\partial u(x,y,t)}{\partial t},&\gamma =1.
 \end{array}
 \right.
\end{equation}

For the 2nd order discretization of the Riemann-Liouville fractional derivatives (\ref{1.4}) and  (\ref{1.5}), it has been detailedly discussed in \cite{Chen:11} being the sequel of \cite{Sousa:12}. This paper applies the discretization scheme to Riesz fractional derivative. The implicit and explicit finite difference schemes are designed. We theoretically prove that the implicit finite difference scheme is unconditionally stable and the stability condition of the explicit scheme is $\frac{\tau^{\gamma}}{(\Delta x)^{\alpha}} +\frac{\tau^{\gamma}}{(\Delta y)^{\beta}} <C$, confirming the general conclusions on the stability conditions of the explicit schemes for fractional PDEs \cite{Deng:12}. The desired 2nd order convergence in space and $(2-\gamma)$-th order convergence in time of both implicit and explicit schemes are theoretically proved and numerically verified.

The  outline of this paper is   as follows.
In Section 2, we introduce the approximation of the  Caputo  fractional derivative and  the
 2nd order finite difference discretizations for the
Riesz fractional derivatives, and derive the full discretization schemes   of  (\ref{1.1}).
 In Sections 3 and 4, the
stability and convergence of the provided implicit and explicit finite difference schemes are
analyzed, respectively.  To show the effectiveness of the schemes, we
perform the numerical experiments to verify the theoretical results in Section 5.
Finally, we conclude the paper with some remarks in the last section.

\section{Derivation of the finite difference scheme}\label{sec:1}
We use two subsections to derive the full discretization schemes of (\ref{1.1}).
The first subsection introduces the approximation of the  Caputo  fractional derivative and  the
 second order finite difference discretizations for the
Riesz fractional    derivatives  in a finite domain.
The second subsection  gives the full discretization scheme (implicit scheme and explicit schemes)
to  the one-dimensional case of (\ref{1.1})
and  (\ref{1.1}) itself, respectively.

\subsection{Discretizations for the Caputo and Riesz fractional derivatives}
Take the mesh points $x_i=x_L+i\Delta x,i=0,1,\ldots ,{N_x}$, $y_j=y_L+j\Delta y,j=0,1,\ldots ,{N_y}$
and $t_k=k\tau,k=0,1,\ldots ,{N_t}$, where
 $\Delta x=(x_R-x_L)/{N_x}$, $\Delta y=(y_R-y_L)/{N_y}$, $\tau=T/{N_t}$,
 i.e., $\Delta x$ and $\Delta y$ are the uniform space stepsizes in the corresponding directions, $\tau$ the time stepsize.
 For $\nu\in (1,2)$,  the left and right  Riemann-Liouville space fractional
derivatives (\ref{1.4}) and (\ref{1.5}) have the second order approximation operators
$\delta_{\nu,_+x}{u_{i,j}^k}$ and $\delta_{\nu,_-x}{u_{i,j}^k}$, respectively, given in a finite
domain \cite{Chen:11}, where ${u_{i,j}^k}$ denotes the
approximated value of $u(x_i,y_j,t_k)$.

The approximation operator of (\ref{1.4}) is defined by \cite{Chen:11}
\begin{equation}\label{2.7}
  \delta_{\nu,_+x}{u_{i,j}^k}:=\frac{1}{\Gamma(4-\nu)(\Delta
  x)^{\nu}}\sum_{m=0}^{i+1}u_{m,j}^kp_{i,m}^{\nu},
\end{equation}
and there exists
\begin{equation}\label{2.8}
  _{x_L}\!D_x^{\nu}u(x,y,t)= \delta_{\nu,_+x}{u_{i,j}^k}+\mathcal{O}(\Delta x)^2,
\end{equation}
where
\begin{equation}\label{2.9}
p_{i,m}^{\nu}=\left\{ \begin{array}
 {l@{\quad } l}
  a_{i-1,m}-2a_{i,m}+a_{i+1,m},&m \leq i-1,\\
  -2a_{i,i}+a_{i+1,i} ,&m=i,\\
 a_{i+1,i+1} , & m=i+1,\\
 0 , &m>i+1,
 \end{array}
 \right.
\end{equation}
and
\begin{equation*}
 a_{j,m}=\left\{ \begin{array}
 {l@{\quad} l}
 (j-1)^{3-\nu}-j^{2-\nu}(j-3+\nu),& m=0,\\
  (j-m+1)^{3-\nu}-2(j-m)^{3-\nu}+(j-m-1)^{3-\nu},&1 \leq m \leq j-1,\\
  1,&m=j,
 \end{array}
 \right.
\end{equation*}
with $j=i-1,i,i+1$.

Analogously,   the   approximation operator of (\ref{1.5}) is  described as \cite{Chen:11}
 \begin{equation}\label{2.10}
  \delta_{\nu,_-x}{u_{i,j}^k}:=\frac{1}{\Gamma(4-\nu)(\Delta
  x)^{\nu}}\sum_{m=i-1}^{N_x}u_{m,j}^kq_{i,m}^{\nu},
\end{equation}
and it holds that
\begin{equation}\label{2.11}
  _{x}D_{x_R}^{\nu}u(x,y,t)= \delta_{\nu,_-x}{u_{i,j}^k}+\mathcal{O}(\Delta x)^2,
\end{equation}
with
\begin{equation}\label{2.12}
q_{i,m}^{\nu}=\left\{ \begin{array}
 {l@{\quad } l}
   0,&m<i-1, \\
  b_{i-1,i-1},&m=i-1,\\
 -2b_{i,i}+b_{i-1,i}, & m=i,\\
 b_{i-1,m}-2b_{i,m}+b_{i+1,m}, &i+1\leq m \leq {N_x},
 \end{array}
 \right.
\end{equation}
and
\begin{equation*}
 b_{j,m}=\left\{ \begin{array}
 {l@{\quad } l}
   1,&m=j ,\\
  (m-j+1)^{3-\alpha}-2(m-j)^{3-\alpha}+(m-j-1)^{3-\alpha},&j+1 \leq m \leq {N_x}-1,\\
 (3-\alpha-{N_x}+j)({N_x}-j)^{2-\alpha}+({N_x}-j-1)^{3-\alpha},& m={N_x},
 \end{array}
 \right.
\end{equation*}
with $j=i-1,i,i+1$.

Combining (\ref{2.8}) and (\ref{2.11}), we obtain the approximation operator of the Riesz fractional derivative

\begin{equation}\label{2.13}
\begin{split}
\frac{\partial ^{\nu}u(x,y_j,t_k)}{\partial |x|^{\nu}}\Big|_{x=x_i}
 &=-\kappa_{\nu}\left( _{x_L}\!D_x^{\nu}+ _{x}\!D_{x_R}^{\nu} \right)u(x,y_j,t_k)\big|_{x=x_i}\\
& =
 -\kappa_{\nu} \left( \delta_{\nu,_+x}+\delta_{\nu,_-x} \right) {u_{i,j}^k} +\mathcal{O}(\Delta x)^2\\
&=\frac{-\kappa_{\nu}}{\Gamma(4-\nu)\Delta
  x^{\nu}}\sum_{m=0}^{N_x} \left( p_{i,m}^{\nu}+q_{i,m}^{\nu}\right)u_{m,j}^k+\mathcal{O}(\Delta x)^2,\\
  &:=\frac{-\kappa_{\nu}}{\Gamma(4-\nu)\Delta
  x^{\nu}}\sum_{m=0}^{N_x}  g_{i,m}^{\nu}u_{m,j}^k+\mathcal{O}(\Delta x)^2,
\end{split}
\end{equation}
where
\begin{equation}\label{2.14}
g_{i,m}^{\nu}=\left\{ \begin{array}
 {l@{\quad } l}
  p_{i,m}^{\nu},&m < i-1\\
  p_{i,i-1}^{\nu}+q_{i,i-1}^{\nu} ,&m=i-1\\
 p_{i,i}^{\nu}+q_{i,i}^{\nu} ,&m=i\\
 p_{i,i+1}^{\nu}+q_{i,i+1}^{\nu} ,&m=i+1\\
 q_{i,m}^{\nu} ,&m>i+1.\\
 \end{array}
 \right.
\end{equation}
 Taking  $\nu=2$, both  Eq. (\ref{2.8}) and (\ref{2.11}) reduce to the following form
\begin{equation*}
 \frac{\partial^2 u(x_i,y,t)}{\partial x^2}=\frac{u(x_{i+1},y,t)-2u(x_i,y,t)+u(x_{i-1},y,t)}{(\Delta x)^2}+\mathcal{O}(\Delta x)^2.
\end{equation*}

The Caputo derivative in the time direction is discretized as  \cite{Lin:07}
\begin{equation}\label{2.15}
  \begin{split}
    _0^C\!D_t^{\gamma}u(x_i,y_j,t)\big|_{t=t_{k+1}}=&\frac{1}{\Gamma(1-\gamma)}\sum_{s=0}^{k}\int_{t_s}^{t_{s+1}}
    \frac{\partial u(x_i,y_j,\eta)}{\partial \eta} \frac{d\eta}{(t_{k+1}-\eta)^{\gamma}}\\
                               =&\frac{1}{\Gamma(1-\gamma)}\sum_{s=0}^{k}\frac{u(x_i,y_j,t_{s+1})-u(x_i,y_j,t_{s})}{\tau}
                               \int_{t_s}^{t_{s+1}}\frac{d\eta}{(t_{k+1}-\eta)^{\gamma}}+\mathcal{O}(\tau^{2-\gamma})\\
                               =&\frac{1}{\Gamma(2-\gamma)}\sum_{s=0}^{k}l_s\frac{u(x_i,y_j,t_{k+1-s})-u(x_i,y_j,t_{k-s})}{\tau^{\gamma}}
                               +\mathcal{O}(\tau^{2-\gamma}),\\
  \end{split}
\end{equation}
where $\gamma \in (0,1)$, $l_s=(s+1)^{1-\gamma}-s^{1-\gamma}$.

When $0<\gamma<1$,  the time Caputo fractional derivative  uses the information of the classical derivatives
 at all previous time levels (non-Markovian process).  If  $\gamma=1$, then $l_0=1$, $l_s=0,s>0$, it can be seen that by taking the limit $\gamma \rightarrow 1$ in (\ref{2.15}), which gives the following equation
\begin{equation*}
\frac{\partial u(x,y,t_{k+1})}{\partial t} =\frac{u(x,y,t_{k+1})-u(x,y,t_{k})}{\tau} + \mathcal{O}(\tau).
\end{equation*}

Similarly, it is easy to get the one-dimensional case of (\ref{2.7})-(\ref{2.15}).

\noindent{\bf Remark 2.1 }\cite{Chen:11}.  Denoting $\tilde{U}^n=[u_{1,j}^n,u_{2,j}^n,\cdots,u_{N_x-1,j}^n]^{\rm T}$,
$j=0,1,\ldots ,{N_y}$ and
 rewriting (\ref{2.7}) and (\ref{2.10}) as matrix forms
$\delta_{\alpha,_+x} \tilde{U}^n=\tilde{A}\tilde{U}^n+b_1$ and
$\delta_{\alpha,_-x} \tilde{U}^n=\tilde{B}\tilde{U}^n+b_2$,
respectively, then there exists $\tilde{A}=\tilde{B}^{\rm T}$.

\subsection{Implicit and explicit difference schemes}

Let us first consider  the one-dimensional time-space Caputo-Riesz fractional   diffusion equation
\begin{equation}\label{2.16}
 _0^C\!D_t^{\gamma}u(x,t)=c(x,t)\frac{\partial ^{\alpha}u(x,t)}{\partial |x|^{\alpha}}+f(x,t).
\end{equation}
 Using the one-dimensional case of (\ref{2.7})-(\ref{2.15}), we can write (\ref{2.16}) as
\begin{equation*}
  \begin{split}
&\frac{{\tau^{-\gamma}}}{\Gamma(2-\gamma)}\sum_{s=0}^{k}l_s [u(x_i,t_{k+1-s})-u(x_i,t_{k-s})]\\
&\quad =-\frac{\kappa_{\alpha}c(x_i,t_{k+1})}{\Gamma(4-\alpha)(\Delta x)^{\alpha}}\sum_{m=0}^{N_x}  g_{i,m}^{\alpha}u(x_m,t_{k+1})+f(x_i,t_{k+1})+\mathcal{O} (\tau^{2-\gamma}+(\Delta x)^2),
  \end{split}
\end{equation*}
where  $i=0,1,\ldots ,{N_x}$, $k=0,1,\ldots ,{N_t}$.
Assuming that $c_i^k=c(x_i,t_k)$, $f_i^k=f(x_i,t_k)$,  $\omega_{i,k+1}=-\frac{\Gamma(2-\gamma) \tau^{\gamma}\kappa_{\alpha}c_i^{k+1}}{\Gamma(4-\alpha)(\Delta x)^{\alpha}}$ and $\mu=\Gamma(2-\gamma)\tau^{\gamma}$,
 we have
\begin{equation}\label{2.17}
  \begin{split}
u(x_i,t_{k+1})&=u(x_i,t_{k})-\sum_{s=1}^{k}l_s[u(x_i,t_{k+1-s})-u(x_i,t_{k-s})]\\
&\quad +\omega_{i,k+1}\sum_{m=0}^{N_x}  g_{i,m}^{\alpha}u(x_m,t_{k+1})+\mu f(x_i,t_{k+1})+R_i^{k+1},
  \end{split}
\end{equation}
where $|R_i^{k+1}|\leq  C\tau^{\gamma}(\tau^{2-\gamma}+(\Delta x)^2)$.

Therefore,  the implicit  difference scheme of (\ref{2.16}) has the following form
\begin{equation}\label{2.18}
  \begin{split}
u_i^{k+1}=u_i^k-\sum_{s=1}^{k}l_s(u_i^{k+1-s}-u_i^{k-s})+\omega_{i,k+1}\sum_{m=0}^{N_x}  g_{i,m}^{\alpha}u_m^{k+1}+\mu f_i^{k+1},
  \end{split}
\end{equation}
and it can be rewrote as
\begin{equation}\label{2.19}
  \begin{split}
(1-\omega_{i,1} g_{i,i}^{\alpha})u_i^1-\omega_{i,1}\!\!\!\!\sum_{m=0,m\neq i}^{N_x} \!\!\!\!  g_{i,m}^{\alpha}u_m^{1}&=u_i^{0}+\mu f_i^{1}, ~~~k=0,\\
(1-\omega_{i,k+1} g_{i,i}^{\alpha})u_i^{k+1}-\omega_{i,k+1}\!\!\!\!\sum_{m=0,m\neq i}^{N_x} \!\!\!\!  g_{i,m}^{\alpha}u_m^{k+1}
  &=\sum_{s=0}^{k-1}(l_s-l_{s+1})u_i^{k-s}+l_ku_i^{0}+\mu f_i^{k+1},    ~~k>0.
  \end{split}
  \end{equation}

Taking $\gamma=1$, thus the implicit  difference scheme (\ref{2.19}) reduces to the following equation
  \begin{equation*}
(1-\omega_{i,k+1}g_{i,i}^{\alpha})u_i^{k+1}-\omega_{i,k+1}\!\!\!\!\sum_{m=0,m\neq i}^{N_x} \!\!\!\! g_{i,m}^{\alpha}u_m^{k+1}
 =u_i^{k}+\mu f_i^{k+1},\quad k\geq0.
  \end{equation*}

 Next,  we examine   the two-dimensional time-space Caputo-Riesz fractional  diffusion equation (\ref{1.1}). According to  (\ref{2.7})-(\ref{2.15}),  then (\ref{1.1}) can be recast as
\begin{equation*}
  \begin{split}
&\frac{{\tau^{-\gamma}}}{\Gamma(2-\gamma)}\!\sum_{s=0}^{k}l_s [u(x_i,y_j,t_{k+1-s})-\!u(x_i,y_j,t_{k-s})]
 =-\frac{\kappa_{\alpha}c(x_i,y_j,t_{k+1})}{\Gamma(4-\alpha)(\Delta x)^{\alpha}}\!\sum_{m=0}^{N_x}  g_{i,m}^{\alpha}u(x_m,y_j,t_{k+1})\\
 &\quad -\frac{\kappa_{\beta}d(x_i,y_j,t_{k+1})}{\Gamma(4-\beta)(\Delta y)^{\beta}}\sum_{m=0}^{N_y}  g_{j,m}^{\beta}u(x_i,y_m,t_{k+1})+f(x_i,y_j,t_{k+1})+\mathcal{O} (\tau^{2-\gamma}+(\Delta x)^2+(\Delta y)^2),
  \end{split}
\end{equation*}
where  $i=0,1,\ldots ,{N_x},j=0,1,\ldots ,{N_y}$, $k=0,1,\ldots ,{N_t}$.
Denoting $c_{i,j}^k=c(x_i,y_j,t_k)$, $d_{i,j}^k=d(x_i,y_j,t_k)$, $f_{i,j}^k=f(x_i,y_j,t_k)$ and
$$
\omega'_{i,j,k+1}=-\frac{\Gamma(2-\gamma) \tau^{\gamma}\kappa_{\alpha}c_{i,j}^{k+1}}{\Gamma(4-\alpha)(\Delta x)^{\alpha}};~~~\omega''_{i,j,k+1}=-\frac{\Gamma(2-\gamma) \tau^{\gamma}\kappa_{\beta}d_{i,j}^{k+1}}{\Gamma(4-\beta)(\Delta y)^{\beta}};
$$
and $\mu=\Gamma(2-\gamma)\tau^{\gamma}$, we obtain
\begin{equation}\label{2.20}
  \begin{split}
u(x_i,y_j,t_{k+1})&=u(x_i,y_j,t_{k})-\sum_{s=1}^{k}l_s[u(x_i,y_j,t_{k+1-s})-u(x_i,y_j,t_{k-s})]\\
 &\quad +\omega'_{i,j,k+1}\sum_{m=0}^{N_x}  g_{i,m}^{\alpha}u(x_m,y_j,t_{k+1})+\omega''_{i,j,k+1}\sum_{m=0}^{N_y}  g_{j,m}^{\beta}u(x_i,y_m,t_{k+1})\\
 &\quad +\mu f(x_i,y_j,t_{k+1})+R_{i,j}^{k+1},
  \end{split}
\end{equation}
where $|R_{i,j}^{k+1}|\leq  C\tau^{\gamma}(\tau^{2-\gamma}+(\Delta x)^2+(\Delta y)^2)$.

Then we obtain the full discretization implicit difference scheme of (\ref{1.1}) as
\begin{equation}\label{2.21}
u_{i,j}^{k+1}=u_{i,j}^k-\sum_{s=1}^{k}l_s(u_{i,j}^{k+1-s}-u_{i,j}^{k-s})+\omega'_{i,j,k+1}\sum_{m=0}^{N_x}  g_{i,m}^{\alpha}u_{m,j}^{k+1}
+\omega''_{i,j,k+1}\sum_{m=0}^{N_y}  g_{j,m}^{\beta}u_{i,m}^{k+1}+\mu f_{i,j}^{k+1},
 \end{equation}
and Eq. (\ref{2.21}) can be rewritten as
\begin{equation}\label{2.22}
  \begin{split}
&(1-\omega'_{i,j,1} g_{i,i}^{\alpha}-\omega''_{i,j,1} g_{j,j}^{\beta})u_{i,j}^{1}-\omega'_{i,j,1}\!\!\!\!\sum_{m=0,m\neq i}^{N_x}  g_{i,m}^{\alpha}u_{m,j}^{1}
-\omega''_{i,j,1}\!\!\!\!\sum_{m=0,m\neq j}^{N_y}  g_{j,m}^{\beta}u_{i,m}^{1}\\
&\qquad\qquad\qquad\qquad=u_{i,j}^{0}+\mu f_{i,j}^{1},   \quad k=0,\\
&(1-\omega'_{i,j,k+1} g_{i,i}^{\alpha}-\omega''_{i,j,k+1} g_{j,j}^{\beta})u_{i,j}^{k+1}-\omega'_{i,j,k+1}\!\!\!\!\sum_{m=0,m\neq i}^{N_x}  g_{i,m}^{\alpha}u_{m,j}^{k+1}
-\omega''_{i,j,k+1}\!\!\!\!\sum_{m=0,m\neq j}^{N_y}  g_{j,m}^{\beta}u_{i,m}^{k+1}\\
 &\qquad\qquad\qquad\qquad=\sum_{s=0}^{k-1}(l_s-l_{s+1})u_{i,j}^{k-s}+l_ku_{i,j}^{0}+\mu f_{i,j}^{k+1},   \quad k>0.
  \end{split}
 \end{equation}
When $\gamma=1$, (\ref{2.22}) becomes
\begin{equation*}
\begin{split}
&(1-\omega'_{i,j,k+1} g_{i,i}^{\alpha}-\omega''_{i,j,k+1} g_{j,j}^{\beta})u_{i,j}^{k+1}-\omega'_{i,j,k+1}\!\!\!\!\sum_{m=0,m\neq i}^{N_x}  \!\!\!\!g_{i,m}^{\alpha}u_{m,j}^{k+1}
-\omega''_{i,j,k+1}\!\!\!\!\sum_{m=0,m\neq j}^{N_y} \!\!\!\! g_{j,m}^{\beta}u_{i,m}^{k+1}\\
&\qquad\qquad\qquad\qquad=u_{i,j}^{k}+\mu f_{i,j}^{k+1},   \quad k\geq0.
\end{split}
\end{equation*}

Let $ \mathbf{u}^k=[\mathbf{u}_0^k,\mathbf{u}_1^k,\ldots,\mathbf{u}_{N_x}^k]^T$,
$ \mathbf{f}^k=[\mathbf{f}_0^k,\mathbf{f}_1^k,\ldots,\mathbf{f}_{N_x}^k]^T$,
where $ \mathbf{u}_i^k=[u_{i,0}^k,u_{i,1}^k,\ldots,u_{i,{N_y}}^k]^T$
and $ \mathbf{f}_i^k=[f_{i,0}^k,f_{i,1}^k,\ldots,f_{i,{N_y}}^k]^T$.
Then (\ref{2.22}) can be written in the matrix form
\begin{equation}\label{2.23}
  \begin{split}
A\mathbf{u}^1&=\mathbf{u}^0+\mu \mathbf{f}^1,   \quad k=0,\\
A\mathbf{u}^{k+1}&=\sum_{s=0}^{k-1}(l_s-l_{s+1})\mathbf{u}^{k-s}+l_k\mathbf{u}^0+\mu \mathbf{f}^{k+1},  \quad  k>0,
  \end{split}
  \end{equation}
where $A $ is a $(N_xN_y)\times (N_xN_y)$  coefficient  matrix.

Analogously,  the explicit  difference scheme of (\ref{2.16}) is
  \begin{equation}\label{2.24}
  \begin{split}
&u_i^1=(1+\sigma_{i,0} g_{i,i}^{\alpha})u_i^{0}+\sigma_{i,0}\!\!\!\!\sum_{m=0,m\neq i}^{N_x} \!\!\!\! g_{i,m}^{\alpha}u_m^{0}+\mu f_i^{0}, \qquad k=0,\\
&u_i^{k+1}=(1-l_1+\sigma_{i,k} g_{i,i}^{\alpha})u_i^{k}+\sum_{s=1}^{k-1}(l_s-l_{s+1})u_i^{k-s}+l_ku_i^{0}
+\sigma_{i,k}\!\!\!\!\sum_{m=0,m\neq i}^{N_x} \!\!\!\! g_{i,m}^{\alpha}u_m^{k}+\mu f_i^{k},    ~~~ k>0,
\end{split}
  \end{equation}
where $\sigma_{i,k}=-\frac{\Gamma(2-\gamma) \tau^{\gamma}\kappa_{\alpha}c_i^{k}}{\Gamma(4-\alpha)(\Delta x)^{\alpha}}>0$.
And the explicit difference scheme of (\ref{1.1}) can be expressed as
  \begin{equation}\label{2.25}
  \begin{split}
u_{i,j}^{1}&=(1+\sigma'_{i,j,0} g_{i,i}^{\alpha}+\sigma''_{i,j,0} g_{j,j}^{\beta})u_{i,j}^{0}
+\sigma'_{i,j,0}\!\!\!\!\sum_{m=0,m\neq i}^{N_x}  g_{i,m}^{\alpha}u_{m,j}^{0}\\
&\quad+\sigma''_{i,j,0}\!\!\!\!\sum_{m=0,m\neq j}^{N_y}  g_{j,m}^{\beta}u_{i,m}^{0}+\mu f_{i,j}^{0}, \qquad ~~~~~~~~~~~~~~~~~~~~~~~~~~~~~~~~~~ k=0,\\
u_{i,j}^{k+1}&=(1-l_1+\sigma'_{i,j,k} g_{i,i}^{\alpha}+\sigma''_{i,j,k} g_{j,j}^{\beta})u_{i,j}^{k}+\sum_{s=1}^{k-1}(l_s-l_{s+1})u_{i,j}^{k-s}+l_ku_{i,j}^{0}\\
&\quad+\sigma'_{i,j,k}\!\!\!\!\sum_{m=0,m\neq i}^{N_x}  g_{i,m}^{\alpha}u_{m,j}^{k}
+\sigma''_{i,j,k}\!\!\!\!\sum_{m=0,m\neq j}^{N_y}  g_{j,m}^{\beta}u_{i,m}^{k}+\mu f_{i,j}^{k},   \qquad ~~~~~~ k>0,
\end{split}
  \end{equation}
where $\sigma'_{i,j,k}=-\frac{\Gamma(2-\gamma) \tau^{\gamma}\kappa_{\alpha}c_{i,j}^{k}}{\Gamma(4-\alpha)(\Delta x)^{\alpha}}$,
 $\sigma''_{i,j,k}=-\frac{\Gamma(2-\gamma) \tau^{\gamma}\kappa_{\beta}d_{i,j}^{k}}{\Gamma(4-\beta)(\Delta y)^{\beta}}$.

\section{Stability analysis}
Now we perform the detailed stability analysis for the implicit and explicit schemes (\ref{2.19}) and (\ref{2.24}) of the one dimensional case (\ref{2.16}), and the implicit and explicit schemes (\ref{2.22}) and (\ref{2.25}) of the two dimensional case (\ref{1.1}). First we introduce two lemmas on the properties of the coefficients of the discretized fractional operators.

\noindent{\bf Lemma 3.1 }\cite{Lin:07,Liu:07}.  \label{La:2}
\emph{ Let $\gamma \in (0,1)$, then  coefficients $l_s$ defined in (\ref{2.15}) satisfy}

  \begin{equation*}\label{3.26}
  \begin{split}
&(1) ~~ l_s>0, s=0,1,\ldots, k.\\
&(2) ~~ 1=l_0>l_1 \cdots > l_k, \quad  l_k \rightarrow 0 \quad as \quad k \rightarrow  \infty.\\
&(3) ~~ C_1k^{\gamma}\leq (l_k)^{-1} \leq C_2k^{\gamma}, \quad  \mbox {where $ C_1$ and $C_2$ are constants. }\\
&(4) ~~ \sum_{s=0}^{k}(l_s-l_{s+1})+ l_{k+1}=(1-l_1)+\sum_{s=1}^{k-1}(l_s-l_{s+1})+l_k=1.
\end{split}
  \end{equation*}

\noindent{\bf Lemma 3.2.}\label{Lemma:3.2.}
\emph{The coefficients $g_{i,m}^{\nu}$, $\nu \in(1,2]$ defined in (\ref{2.14}) satisfy}

\begin{equation*}\label{3.27}
  \begin{split}
&(1) ~~ g_{i,i}^{\nu}<0,\,\,\,\,\,\, g_{i,m}^{\nu} > 0\,\,(m\neq i);\\
&(2) ~~ \sum\limits_{m=0}^{N_x}g_{i,m}^{\nu}< 0 ~~ \mbox {and} ~~  -g_{i,i}^{\nu}>\!\!\!\!\!\sum\limits_{m=0,m\neq i}^{N_x}\!\!\!\!g_{i,m}^{\nu}.
  \end{split}
\end{equation*}

\begin{proof}
According to Theorem B in the Appendix of \cite{Chen:11}, and formulas  (\ref{2.9}) and (\ref{2.12}), we obtain
$q_{i,m}^{\nu}> 0$ for $m > i+1$, $p_{i,m}^{\nu}> 0$ for $m < i-1$, and $p_{i,i+1}^{\nu} =q_{i,i-1}^{\nu} =1$,
   $p_{i,i}^{\nu}= q_{i,i}^{\nu}=-4+2^{3-\nu}$. From the definition of $g_{i,m}^{\nu}$ (\ref{2.14}), it is easy to verify that  $g_{i,i+1}^{\nu}=g_{i,i-1}^{\nu}=7-2^{5-\nu}+3^{3-\nu}>0$, $g_{i,i}^{\nu}=-8+2^{4-\nu}<0$, and other $g_{i,m}^{\nu}$ are positive. And
letting $u$ be a constant in (\ref{2.13}), we can easily derive that $\sum\limits_{m=0}^{N_x}g_{i,m}^{\nu}< 0$, then $-g_{i,i}^{\nu}>\!\!\!\!\!\sum\limits_{m=0,m\neq i}^{N_x}\!\!\!\!g_{i,m}^{\nu}$.
\end{proof}
Next, we use four subsections to strictly prove that both the implicit schemes (\ref{2.19}) and (\ref{2.22}) are unconditionally stable; the explicit scheme (\ref{2.24}) is stable under the condition $\frac{\tau^{\gamma}}{(\Delta x)^{\alpha}}<C$  and the explicit scheme (\ref{2.25}) is stable under the condition $\frac{\tau^{\gamma}}{(\Delta x)^{\alpha}} +\frac{\tau^{\gamma}}{(\Delta y)^{\beta}} <C$.

\subsection{Stability for one-dimensional implicit difference scheme}

\noindent{\bf Theorem 3.3.}\label{Theorem:5}
\emph{The  implicit difference scheme (\ref{2.19}) of the one-dimensional  time-space Caputo-Riesz fractional diffusion equation (\ref{2.16})
with  $0<\gamma \leq1$, $1<\alpha \leq2$ is unconditionally stable.}

\begin{proof}
Let $\widetilde{u_i^k}~(i=0,1,\ldots,N_x;\,k=0,1,\ldots,N_t)$ be the approximate solution of $u_i^k$,
which is the exact solution of the  implicit scheme (\ref{2.19}).
Putting $\epsilon_i^k=\widetilde{u_i^k}-u_i^k$, then from (\ref{2.19}) we obtain the following perturbation equation
\begin{equation}\label{3.28}
  \begin{split}
(1-\omega_{i,1}  g_{i,i}^{\alpha})\epsilon_i^1-\omega_{i,1} \!\!\!\!\sum_{m=0,m\neq i}^{N_x}\!\!\!\!  g_{i,m}^{\alpha}\epsilon_m^1&=\epsilon_i^{0}, \quad k=0,\\
(1-\omega_{i,k+1}g_{i,i}^{\alpha})\epsilon_i^{k+1}-\omega_{i,k+1}\!\!\!\!\!\sum_{m=0,m\neq i}^{N_x} \!\!\!\!  g_{i,m}^{\alpha}\epsilon_m^{k+1}
  &=\sum_{s=0}^{k-1}(l_s-l_{s+1})\epsilon_i^{k-s}+l_k\epsilon_i^{0},    \quad k>0.
  \end{split}
  \end{equation}

When $\gamma=1$, Eq. (\ref{3.28}) can be written as
\begin{align*}
(1-\omega_{i,k+1} g_{i,i}^{\alpha})\epsilon_i^{k+1}-\omega_{i,k+1}\!\!\!\!\!\sum_{m=0,m\neq i}^{N_x} \!\!\!\! g_{i,m}^{\alpha}\epsilon_m^{k+1}
  =\epsilon_i^{k},    \quad k \geq 0.
\end{align*}
Denoting $E^k=[\epsilon_0^k,\epsilon_1^k,\ldots, \epsilon_{N_x}^k]$
and $||E^k||_{\infty}=\max \limits _{0\leq i \leq N_x}|\epsilon_i^k|$, then we use the mathematical induction to prove the unconditional stability.
For  $k=0$, supposing  $|\epsilon_{i_0}^1|=||E^1||_{\infty}=\max \limits _{0\leq i \leq N_x}|\epsilon_i^1|$,
according to Lemma 3.2, we get
 \begin{align*}
||E^1||_{\infty}&=|\epsilon_{i_0}^1|  \leq |\epsilon_{i_0}^1| -\omega_{i_0,1}\sum_{m=0}^{N_x}  g_{{i_0},m}^{\alpha}|\epsilon_{i_0}^{1}|
                     =|\epsilon_{i_0}^1|-\omega_{i_0,1} g_{{i_0},{i_0}}^{\alpha}|\epsilon_{i_0}^{1}|
                                     -\omega_{i_0,1}\!\!\!\!\!\sum_{m=0,m\neq {i_0}}^{N_x} \!\!\!\! g_{{i_0},m}^{\alpha}|\epsilon_{i_0}^{1}|\\
                 &\leq (1-\omega_{i_0,1} g_{{i_0},{i_0}}^{\alpha})|\epsilon_{i_0}^{1}|
                                    -\omega_{i_0,1}\!\!\!\!\!\!\sum_{m=0,m\neq {i_0}}^{N_x} \!\!\!\!\! g_{{i_0},m}^{\alpha}|\epsilon_m^{1}|
                    \leq \big|(1-\omega_{i_0,1} g_{{i_0},{i_0}}^{\alpha})\epsilon_{i_0}^{1}
                                    -\omega_{i_0,1}\!\!\!\!\!\sum_{m=0,m\neq {i_0}}^{N_x}\!\!\!\!  g_{{i_0},m}^{\alpha} \epsilon_m^{1}\big|\\
                 & =|\epsilon_{i_0}^{0}|  \leq ||E^0||_{\infty}.
  \end{align*}
Assuming $||E^{\widetilde k}||_{\infty} \leq ||E^0||_{\infty},\,\, {\widetilde k}=1,2,\ldots,k,$
and $|\epsilon_{i_0}^{k+1}|=||E^{k+1}||_{\infty}=\max \limits _{0\leq i \leq N_x}|\epsilon_i^{k+1}|$,
there exists
  \begin{align*}
||E^{k+1}||_{\infty}&=|\epsilon_{i_0}^{k+1}|  \leq |\epsilon_{i_0}^{k+1}| -\omega_{i_0,k+1}\sum_{m=0}^{N_x}  g_{{i_0},m}^{\alpha}|\epsilon_{i_0}^{k+1}|\\
                    &\quad  =|\epsilon_{i_0}^{k+1}|- \omega_{i_0,k+1} g_{{i_0},{i_0}}^{\alpha}|\epsilon_{i_0}^{k+1}|
                      -\omega_{i_0,k+1}\!\!\!\!\!\sum_{m=0,m\neq {i_0}}^{N_x} \!\!\!\!\! g_{{i_0},m}^{\alpha}|\epsilon_{i_0}^{k+1}|\\
                    &\quad \leq (1-\omega_{i_0,k+1} g_{{i_0},{i_0}}^{\alpha})|\epsilon_{i_0}^{k+1}| -\omega_{i_0,k+1}\!\!\!\!\!\sum_{m=0,m\neq {i_0}}^{N_x}  \!\!\!\!\! g_{{i_0},m}^{\alpha}|\epsilon_m^{k+1}| \\
                     &\quad  \leq \big|(1-\omega_{i_0,k+1} g_{{i_0},{i_0}}^{\alpha})\epsilon_{i_0}^{k+1} -\omega_{i_0,k+1}\!\!\!\!\!\sum_{m=0,m\neq {i_0}}^{N_x}  \!\!\!\!\! g_{{i_0},m}^{\alpha} \epsilon_m^{k+1}\big|\\
                     &\quad=\big|\sum_{s=0}^{k-1}(l_s-l_{s+1})\epsilon_{i_0}^{k-s}+l_k\epsilon_{i_0}^{0}\big|\\
                    &\quad=\big|(1-l_1)\epsilon_{i_0}^{k}+l_k\epsilon_{i_0}^{0} +\sum_{s=1}^{k-1}(l_s-l_{s+1})\epsilon_{i_0}^{k-s}\big|.
\end{align*}
From Lemma 3.1, we obtain
  \begin{align*}
||E^{k+1}||_{\infty} &\leq  (1-l_1) ||E^{k}||_{\infty}  +l_k||E^{0}||_{\infty} +\sum_{s=1}^{k-1}(l_s-l_{s+1})||E^{k-s}||_{\infty}\\
                     &\leq  (1-l_1) ||E^{0}||_{\infty}  +l_k||E^{0}||_{\infty} +\sum_{s=1}^{k-1}(l_s-l_{s+1})||E^{0}||_{\infty}\\
                     & \leq ||E^0||_{\infty}.
\end{align*}
When $\gamma=1$, using similar idea, we can prove
\begin{align*}
||E^{k+1}||_{\infty} \leq  |(1-\omega_{i_0,k+1} g_{i_0,i_0}^{\alpha})\epsilon_{i_0}^{k+1}
                  -\omega_{i_0,k+1}\!\!\!\!\!\sum_{m=0,m\neq i_0}^{N_x} \!\!\!\! g_{i_0,m}^{\alpha}\epsilon_m^{k+1}|
  =|\epsilon_{i_0}^{k}|\leq ||E^0||_{\infty}.
\end{align*}
\end{proof}

\subsection{Stability for one-dimensional explicit difference scheme}

\noindent{\bf Theorem 3.4.}\label{Theorem:07}
\emph{If
\begin{equation*}
0<\frac{\tau^{\gamma}}{(\Delta x)^{\alpha}} \leq -\frac{\Gamma(4-\alpha)(1-2^{-\gamma})}{4\kappa_{\alpha}C_{max}\Gamma(2-\gamma) (1-2^{1-\alpha})},
\quad where \quad  C_{max}=\max_{0\leq i \leq N_x,0\leq k \leq N_t} c(x_i,t_k),
\end{equation*}
 then the explicit difference scheme (\ref{2.24}) of the one-dimensional  time-space Caputo-Riesz fractional   diffusion equation (\ref{2.16})
with  $0<\gamma \leq1$, $1<\alpha \leq2$ is  stable.}

\begin{proof}
Under the above conditions,  we obtain $0<-\sigma_{i,k}g_{i,i}^{\alpha} \leq 2-2^{1-\gamma} $
and  $1-l_1+\sigma_{i,k}g_{i,i}^{\alpha}\geq0$.
Assuming that $\widetilde{u_i^k}~(i=0,1,\ldots,N_x;\,k=0,1,\ldots,N_t)$ be the approximate solution of $u_i^k$,
which is the exact solution of the  explicit scheme (\ref{2.24}).
Therefore the error $\epsilon_i^k=\widetilde{u_i^k}-u_i^k$ satisfies

\begin{equation}\label{3.29}
  \begin{split}
\epsilon_i^1&=(1+\sigma_{i,0} g_{i,i}^{\alpha})\epsilon_i^{0}+\sigma_{i,0}\!\!\!\!\!\sum_{m=0,m\neq i}^{N_x}\!\!\!\!\!  g_{i,m}^{\alpha}\epsilon_m^{0}, \qquad k=0,\\
\epsilon_i^{k+1}&=(1-l_1+\sigma_{i,k} g_{i,i}^{\alpha})\epsilon_i^{k}+\sum_{s=1}^{k-1}(l_s-l_{s+1})\epsilon_i^{k-s}+l_k\epsilon_i^{0}
+\sigma_{i,k}\!\!\!\!\!\sum_{m=0,m\neq i}^{N_x}\!\!\!\!\!  g_{i,m}^{\alpha}\epsilon_m^{k},   \qquad  k>0.
 \end{split}
\end{equation}
When $\gamma=1$, (\ref{3.29}) becomes
\begin{align*}
\epsilon_i^{k+1}=(1+\sigma_{i,k} g_{i,i}^{\alpha})\epsilon_i^{k}
+\sigma_{i,k}\!\!\!\!\!\sum_{m=0,m\neq i}^{N_x} \!\!\!\!\!  g_{i,m}^{\alpha}\epsilon_m^{k},   \quad k \geq 0.
\end{align*}
Denoting $E^k=[\epsilon_0^k,\epsilon_1^k,\ldots, \epsilon_{N_x}^k]$
and $||E^k||_{\infty}=\max \limits _{0\leq i \leq N_x}|\epsilon_i^k|$, we use mathematical induction to prove the conditional stability.
For  $k=0$, supposing  $|\epsilon_{i_0}^1|=||E^1||_{\infty}=\max \limits _{0\leq i \leq N_x}|\epsilon_i^1|$,
according to Lemma 3.2, we get
  \begin{align*}
~~~~~~~~~||E^1||_{\infty}=|\epsilon_{i_0}^1|  &=\big|(1+\sigma_{i_0,0} g_{{i_0},{i_0}}^{\alpha})\epsilon_{i_0}^{0}
                      +\sigma_{i_0,0}\!\!\!\!\!\sum_{m=0,m\neq {i_0}}^{N_x} \!\!\!\!\! g_{{i_0},m}^{\alpha}\epsilon_m^{0}\big|\\
& \leq  (1+\sigma_{i_0,0} g_{{i_0},{i_0}}^{\alpha})|\epsilon_{i_0}^{0}|
                      +\sigma_{i_0,0}\!\!\!\!\!\sum_{m=0,m\neq {i_0}}^{N_x} \!\!\!\!\! g_{{i_0},m}^{\alpha}|\epsilon_m^{0}|\\
&\leq  (1+\sigma_{i_0,0} g_{{i_0},{i_0}}^{\alpha})||E^0||_{\infty}
                      +\sigma_{i_0,0}\!\!\!\!\!\sum_{m=0,m\neq {i_0}}^{N_x} \!\!\!\!\! g_{{i_0},m}^{\alpha}||E^0||_{\infty}\\
&=||E^0||_{\infty}  +   \sigma_{i_0,0}\sum_{m=0}^{N_x}  g_{{i_0},m}^{\alpha}||E^0||_{\infty}
 \leq   ||E^0||_{\infty} .
  \end{align*}
Assuming $||E^{\widetilde k}||_{\infty} \leq ||E^0||_{\infty},\, {\widetilde k}=1,2,\ldots,k,$
and denoting  $|\epsilon_{i_0}^{k+1}|=||E^{k+1}||_{\infty}=\max \limits _{0\leq i \leq N_x}|\epsilon_i^{k+1}|$,
from Lemma 3.1, there exists
  \begin{align*}
||E^{k+1}||_{\infty} &=  |\epsilon_{i_0}^{k+1}|  =\Big|(1-l_1+\sigma_{i_0,k} g_{{i_0},{i_0}}^{\alpha})\epsilon_{i_0}^{k}
                            +\sum_{s=1}^{k-1}(l_s-l_{s+1})\epsilon_{i_0}^{k-s}+l_k\epsilon_{i_0}^{0}
                            +\sigma_{i_0,k}\!\!\!\!\!\sum_{m=0,m\neq {i_0}}^{N_x}\!\!\!\!\!  g_{{i_0},m}^{\alpha}\epsilon_m^{k}\Big|\\
                     &\leq (1-l_1+\sigma_{i_0,k} g_{{i_0},{i_0}}^{\alpha})|\epsilon_{i_0}^{k}|
                            +\sum_{s=1}^{k-1}(l_s-l_{s+1})|\epsilon_{i_0}^{k-s}|+l_k|\epsilon_{i_0}^{0}|
                            +\sigma_{i_0,k}\!\!\!\!\!\sum_{m=0,m\neq {i_0}}^{N_x} \!\!\!\!\! g_{{i_0},m}^{\alpha}|\epsilon_m^{k}|\\
                     &\leq (1-l_1)||E^{k}||_{\infty}
                            +\sum_{s=1}^{k-1}(l_s-l_{s+1})||E^{k-s}||_{\infty}+l_k||E^{0}||_{\infty}
                            +\sigma_{i_0,k}\sum_{m=0}^{N_x}  g_{{i_0},m}^{\alpha}||E^{k}||_{\infty}\\
                     &\leq (1-l_1)||E^{0}||_{\infty}
                            +\sum_{s=1}^{k-1}(l_s-l_{s+1})||E^{0}||_{\infty}+l_k||E^{0}||_{\infty}\\
                     &=  (1-l_1+\sum_{s=1}^{k-1}(l_s-l_{s+1})+l_k)||E^{0}||_{\infty} =||E^{0}||_{\infty}.
\end{align*}
If $\gamma=1$, using similar method, we can prove
\begin{align*}
||E^{k+1}||_{\infty} =|\epsilon_{i_0}^{k+1}|&=\Big|(1+\sigma_{i_0,k} g_{{i_0},{i_0}}^{\alpha})\epsilon_{i_0}^{k}
                                          +\sigma_{i_0,k}\!\!\!\!\!\sum_{m=0,m\neq {i_0}}^{N_x}  g_{{i_0},m}^{\alpha}\epsilon_m^{k}\Big|\\
                                        &\leq (1+\sigma_{i_0,k} g_{{i_0},{i_0}}^{\alpha})|\epsilon_{i_0}^{k}|
                                          +\sigma_{i_0,k}\!\!\!\!\!\sum_{m=0,m\neq {i_0}}^{N_x}  g_{{i_0},m}^{\alpha}|\epsilon_m^{k}|\\
                                        &\leq (1+\sigma_{i_0,k} g_{{i_0},{i_0}}^{\alpha})||E^{k}||_{\infty}
                                          +\sigma_{i_0,k}\!\!\!\!\!\sum_{m=0,m\neq {i_0}}^{N_x}  g_{{i_0},m}^{\alpha}||E^{k}||_{\infty}\\
                                        &=||E^{k}||_{\infty}+\sigma_{i_0,k}\sum_{m=0}^{N_x}  g_{{i_0},m}^{\alpha}||E^{k}||_{\infty}
                                          \leq ||E^{k}||_{\infty}
                                          \leq ||E^{0}||_{\infty}.
\end{align*}
\end{proof}

\subsection{Stability for two-dimensional implicit difference scheme}

\noindent{\bf Theorem 3.5.}\label{Theorem:7}
\emph{The  implicit difference scheme (\ref{2.22}) of the two-dimensional time-space Caputo-Riesz fractional   diffusion equation(\ref{1.1})
with  $0<\gamma \leq1$, $1<\alpha, \beta \leq2$ is unconditionally stable.}

\begin{proof}
Let  $\widetilde{u_{i,j}^k}~(i=0,1,\ldots,N_x;j=0,1,\ldots,N_y;k=0,1,\ldots,N_t)$ be the approximate solution of $u_{i,j}^k$,
which is the exact solution of the implicit scheme (\ref{2.22}). Denoting that  $\epsilon_{i,j}^k=\widetilde{u_{i,j}^k}-u_{i,j}^k$,  from (\ref{2.22}) we get the following perturbation equation
\begin{equation}\label{3.30}
  \begin{split}
&(1-\omega'_{i,j,k+1} g_{i,i}^{\alpha}-\omega''_{i,j,k+1} g_{j,j}^{\beta})\epsilon_{i,j}^{1}-\omega'_{i,j,k+1}\!\!\!\!\!\sum_{m=0,m\neq i}^{N_x}  \!\!\!\!\! g_{i,m}^{\alpha}\epsilon_{m,j}^{1}
-\omega''_{i,j,k+1}\!\!\!\!\!\sum_{m=0,m\neq j}^{N_y} \!\!\!\!\! g_{j,m}^{\beta}\epsilon_{i,m}^{1}=\epsilon_{i,j}^{0},  \quad k=0,\\
&(1-\omega'_{i,j,k+1}g_{i,i}^{\alpha}-\omega''_{i,j,k+1} g_{j,j}^{\beta})\epsilon_{i,j}^{k+1}-\omega'_{i,j,k+1}\!\!\!\!\!\sum_{m=0,m\neq i}^{N_x}  \!\!\!\!\! g_{i,m}^{\alpha}\epsilon_{m,j}^{k+1}
-\omega''_{i,j,k+1}\!\!\!\!\!\sum_{m=0,m\neq j}^{N_y}\!\!\!\!\!  g_{j,m}^{\beta}\epsilon_{i,m}^{k+1}\\
 &\qquad\qquad=\sum_{s=0}^{k-1}(l_s-l_{s+1})\epsilon_{i,j}^{k-s}+l_k\epsilon_{i,j}^{0},    \quad k>0.
 \end{split}
\end{equation}
When $\gamma=1$,  (\ref{3.30})  becomes
\begin{equation*}
(1-\omega'_{i,j,k+1} g_{i,i}^{\alpha}-\omega''_{i,j,k+1} g_{j,j}^{\beta})\epsilon_{i,j}^{k+1}-\omega'_{i,j,k+1}\!\!\!\!\!\sum_{m=0,m\neq i}^{N_x}  \!\!\!\!\!g_{i,m}^{\alpha}\epsilon_{m,j}^{k+1}
-\omega''_{i,j,k+1}\!\!\!\!\!\sum_{m=0,m\neq j}^{N_y} \!\!\!\!\! g_{j,m}^{\beta}\epsilon_{i,m}^{k+1}
 =\epsilon_{i,j}^{k},        \quad k \geq 0.
\end{equation*}
Denote $\mathbf{u}^k=[\mathbf{u}_0^k,\mathbf{u}_1^k,\ldots,\mathbf{u}_{N_x}^k]^T$,
 $\mathbf{u}_i^k=[\epsilon_{i,0}^k,\epsilon_{i,1}^k,\ldots,\epsilon_{i,{N_y}}^k]^T$
and $||\mathbf{E}^k||_{\infty}=\max \limits _{{0\leq i \leq N_x},{0\leq j \leq N_y}}|\epsilon_{i,j}^k|$.
We prove the results by mathematical induction. For $k=0$, supposing  $|\epsilon_{i_0,j_0}^1|=||\mathbf{E}^1||_{\infty}=\max \limits _{{0\leq i \leq N_x},{0\leq j \leq N_y}}|\epsilon_{i,j}^1|$, there exists
\begin{align*}
||\mathbf{E}^1||_{\infty}=& |\epsilon_{i_0,j_0}^1|
                     \leq  |\epsilon_{i_0,j_0}^1| -\omega'_{i_0,j_0,1}\sum_{m=0}^{N_x}  g_{{i_0},m}^{\alpha}|\epsilon_{i_0,j_0}^{1}|
                            -\omega''_{i_0,j_0,1}\sum_{m=0}^{N_y}  g_{{j_0},m}^{\beta}|\epsilon_{i_0,j_0}^{1}|\\
                         =&|\epsilon_{i_0,j_0}^1|
                            -\omega'_{i_0,j_0,1} g_{{i_0},{i_0}}^{\alpha}|\epsilon_{i_0,j_0}^{1}|
                            -\omega'_{i_0,j_0,1}\!\!\!\!\!\sum_{m=0,m\neq {i_0}}^{N_x}\!\!\!\!\! g_{{i_0},m}^{\alpha}|\epsilon_{i_0,j_0}^{1}|\\
                           &-\omega''_{i_0,j_0,1} g_{{j_0},{j_0}}^{\beta}|\epsilon_{i_0,j_0}^{1}|
                            -\omega''_{i_0,j_0,1}\!\!\!\!\!\sum_{m=0,m\neq {j_0}}^{N_y}\!\!\!\!\! g_{{j_0},m}^{\beta}|\epsilon_{i_0,j_0}^{1}|\\
                      \leq & (1-\omega'_{i_0,j_0,1} g_{{i_0},{i_0}}^{\alpha}-\omega''_{i_0,j_0,1} g_{{j_0},{j_0}}^{\beta})|\epsilon_{i_0,j_0}^{1}|\\
                           & -\omega'_{i_0,j_0,1}\!\!\!\!\!\sum_{m=0,m\neq {i_0}}^{N_x}\!\!\!\!\!  g_{{i_0},m}^{\alpha}|\epsilon_{m,{j_0}}^{1}|
                            -\omega''_{i_0,j_0,1}\!\!\!\!\!\sum_{m=0,m\neq {j_0}}^{N_y} \!\!\!\!\! g_{{j_0},m}^{\beta}|\epsilon_{{i_0,m}}^{1}|\\
                      \leq & \Big |(1-\omega'_{i_0,j_0,1} g_{{i_0},{i_0}}^{\alpha}-\omega''_{i_0,j_0,1} g_{{j_0},{j_0}}^{\beta})\epsilon_{i_0,j_0}^{1}\\
                           & -\omega'_{i_0,j_0,1}\!\!\!\!\!\sum_{m=0,m\neq {i_0}}^{N_x}  g_{{i_0},m}^{\alpha}\epsilon_{m,{j_0}}^{1}
                            -\omega''_{i_0,j_0,1}\!\!\!\!\!\sum_{m=0,m\neq {j_0}}^{N_y}  g_{{j_0},m}^{\beta}\epsilon_{{i_0,m}}^{1}\Big|\\
                          =&|\epsilon_{i_0,j_0}^{0}| \leq ||\mathbf{E}^0||_{\infty}.
  \end{align*}
Supposing $||\mathbf{E}^{\widetilde k}||_{\infty} \leq ||\mathbf{E}^0||_{\infty},\, {\widetilde k}=1,2,\ldots,k,$
and  $|\epsilon_{i_0,j_0}^{k+1}|=||\mathbf{E}^{k+1}||_{\infty}=\max \limits _{{0\leq i \leq N_x},{0\leq j \leq N_y}}|\epsilon_{i,j}^{k+1}|$,
then
  \begin{align*}
&||\mathbf{E}^{k+1}||_{\infty}= |\epsilon_{i_0,j_0}^{k+1}|  \leq |\epsilon_{i_0,j_0}^{k+1}| -\omega'_{i_0,j_0,k+1}\sum_{m=0}^{N_x}  g_{{i_0},m}^{\alpha}|\epsilon_{i_0,j_0}^{k+1}|
                            -\omega''_{i_0,j_0,k+1}\sum_{m=0}^{N_y}  g_{{j_0},m}^{\beta}|\epsilon_{i_0,j_0}^{k+1}|\\
                         &~=|\epsilon_{i_0,j_0}^{k+1}|
                            -\omega'_{i_0,j_0,k+1}g_{{i_0},{i_0}}^{\alpha}|\epsilon_{i_0,j_0}^{k+1}|
                            -\omega'_{i_0,j_0,k+1}\!\!\!\!\!\sum_{m=0,m\neq {i_0}}^{N_x}\!\!\!\!\! g_{{i_0},m}^{\alpha}|\epsilon_{i_0,j_0}^{k+1}|\\
                          &~\quad  -\omega''_{i_0,j_0,k+1} g_{{j_0},{j_0}}^{\beta}|\epsilon_{i_0,j_0}^{k+1}|
                            -\omega''_{i_0,j_0,k+1}\!\!\!\!\!\sum_{m=0,m\neq {j_0}}^{N_y}\!\!\!\!\! g_{{j_0},m}^{\beta}|\epsilon_{i_0,j_0}^{k+1}|\\
                     &~ \leq  (1-\omega'_{i_0,j_0,k+1} g_{{i_0},{i_0}}^{\alpha}-\omega''_{i_0,j_0,k+1} g_{{j_0},{j_0}}^{\beta})|\epsilon_{i_0,j_0}^{k+1}|\\
                     &~\quad      -\omega'_{i_0,j_0,k+1}\!\!\!\!\!\sum_{m=0,m\neq {i_0}}^{N_x}  g_{{i_0},m}^{\alpha}|\epsilon_{m,{j_0}}^{k+1}|
                            -\omega''_{i_0,j_0,k+1}\!\!\!\!\!\sum_{m=0,m\neq {j_0}}^{N_y}  g_{{j_0},m}^{\beta}|\epsilon_{{i_0,m}}^{k+1}|\\
                     &~ \leq \Big |(1-\omega'_{i_0,j_0,k+1} g_{{i_0},{i_0}}^{\alpha}-\omega''_{i_0,j_0,k+1} g_{{j_0},{j_0}}^{\beta})\epsilon_{i_0,j_0}^{k+1}\\
                     &~\quad -\omega'_{i_0,j_0,k+1}\!\!\!\!\!\sum_{m=0,m\neq {i_0}}^{N_x}  g_{{i_0},m}^{\alpha}\epsilon_{m,{j_0}}^{k+1}
                            -\omega''_{i_0,j_0,k+1}\!\!\!\!\!\sum_{m=0,m\neq {j_0}}^{N_y}  g_{{j_0},m}^{\beta}\epsilon_{{i_0,m}}^{k+1}\Big|\\
                    &~=\big|\sum_{s=0}^{k-1}(l_s-l_{s+1})\epsilon_{i_0,j_0}^{k-s}+l_k\epsilon_{i_0,j_0}^{0}\big|
                      =\big|(1-l_1)\epsilon_{i_0,j_0}^{k}+l_k\epsilon_{i_0,j_0}^{0} +\sum_{s=1}^{k-1}(l_s-l_{s+1})\epsilon_{i_0,j_0}^{k-s}\big|.
\end{align*}
By Lemma 3.1, we get
  \begin{align*}
||\mathbf{E}^{k+1}||_{\infty} &\leq  (1-l_1) ||\mathbf{E}^{k}||_{\infty}  +l_k||\mathbf{E}^{0}||_{\infty}
                                         + \sum_{s=1}^{k-1}(l_s-l_{s+1})||\mathbf{E}^{k-s}||_{\infty}\\
                              &\leq  (1-l_1) ||\mathbf{E}^{0}||_{\infty} +l_k||\mathbf{E}^{0}||_{\infty}
                                         +\sum_{s=1}^{k-1}(l_s-l_{s+1})||\mathbf{E}^{0}||_{\infty}\\
                              &\leq ||\mathbf{E}^0||_{\infty}.
\end{align*}
If $\gamma=1$, we obtain
\begin{align*}
||\mathbf{E}^{k+1}||_{\infty} &\leq \Big|(1-\omega'_{i_0,j_0,k+1} g_{i_0,i_0}^{\alpha}-\omega''_{i_0,j_0,k+1} g_{j_0,j_0}^{\beta})\epsilon_{i_0,j_0}^{k+1}\\
     &~\quad  -\omega'_{i_0,j_0,k+1}\!\!\!\!\!\sum_{m=0,m\neq i}^{N_x}  g_{i_0,m}^{\alpha}\epsilon_{m,j_0}^{k+1}
     -\omega''_{i_0,j_0,k+1}\!\!\!\!\!\sum_{m=0,m\neq j}^{N_y}  g_{j_0,m}^{\beta}\epsilon_{i_0,m}^{k+1}\Big|\\
 &=|\epsilon_{i_0,j_0}^{k}| \leq ||\mathbf{E}^{0}||_{\infty}.
 \end{align*}
\end{proof}

\subsection{Stability for two-dimensional explicit difference scheme}
\noindent{\bf Theorem 3.6.}\label{Theorem:11}
\emph{If
$
0<\frac{\tau^{\gamma}}{(\Delta x)^{\alpha}} +\frac{\tau^{\gamma}}{(\Delta y)^{\beta}}\leq \frac{1-2^{-\gamma}}{\Gamma(2-\gamma)\mathbf{C}_{max} }$,
 where $$\mathbf{C}_{max}=\max_{0\leq i \leq N_x,0\leq j \leq N_y,0\leq k \leq N_t}
\left \{-\frac{\kappa_{\alpha}(4-2^{3-\alpha})c_{i,j}^k}{\Gamma (4-\alpha)},-\frac{\kappa_{\beta}(4-2^{3-\beta})d_{i,j}^k}{\Gamma (4-\beta)}\right\},
$$
then the explicit difference scheme (\ref{2.25}) of the two-dimensional time-space Caputo-Riesz fractional  diffusion equation (\ref{1.1})
with  $0<\gamma \leq1$, $1<\alpha, \beta \leq2$ is stable.}

\begin{proof}
Let  $\widetilde{u_{i,j}^k}~(i=0,1,\ldots,N_x;\,j=0,1,\ldots,N_y;\,k=0,1,\ldots,N_t)$ be the approximate solution of $u_{i,j}^k$,
which is the exact solution of the explicit scheme (\ref{2.25}). Denoting $\epsilon_{i,j}^k=\widetilde{u_{i,j}^k}-u_{i,j}^k$ and using  (\ref{2.25}),  we obtain the following perturbation form
\begin{equation}\label{3.31}
  \begin{split}
\epsilon_{i,j}^{1}&=(1+\sigma'_{i,j,0} g_{i,i}^{\alpha}+\sigma''_{i,j,0} g_{j,j}^{\beta})\epsilon_{i,j}^{0}
+\sigma'_{i,j,0}\!\!\!\!\!\sum_{m=0,m\neq i}^{N_x}  g_{i,m}^{\alpha}\epsilon_{m,j}^{0}+\sigma''_{i,j,0}\!\!\!\!\!\sum_{m=0,m\neq j}^{N_y}  g_{j,m}^{\beta}\epsilon_{i,m}^{0}, \quad k=0,\\
\epsilon_{i,j}^{k+1}&=(1-l_1+\sigma'_{i,j,k} g_{i,i}^{\alpha}+\sigma''_{i,j,k} g_{j,j}^{\beta})\epsilon_{i,j}^{k}+\sum_{s=1}^{k-1}(l_s-l_{s+1})\epsilon_{i,j}^{k-s}+l_k\epsilon_{i,j}^{0}\\
&\quad+\sigma'_{i,j,k}\!\!\!\!\!\sum_{m=0,m\neq i}^{N_x}  g_{i,m}^{\alpha}\epsilon_{m,j}^{k}
       +\sigma''_{i,j,k}\!\!\!\!\!\sum_{m=0,m\neq j}^{N_y}  g_{j,m}^{\beta}\epsilon_{i,m}^{k},    \quad k>0.
  \end{split}
\end{equation}
When $\gamma=1$, (\ref{3.31}) can be rewritten as
\begin{align*}
\epsilon_{i,j}^{k+1}&=(1+\sigma'_{i,j,k} g_{i,i}^{\alpha}+\sigma''_{i,j,k} g_{j,j}^{\beta})\epsilon_{i,j}^{k}
+\sigma'_{i,j,k}\!\!\!\!\!\sum_{m=0,m\neq i}^{N_x}  g_{i,m}^{\alpha}\epsilon_{m,j}^{k}
+\sigma''_{i,j,k}\!\!\!\!\!\sum_{m=0,m\neq j}^{N_y}  g_{j,m}^{\beta}\epsilon_{i,m}^{k}, \quad k \geq 0.
\end{align*}
Using the mathematical induction, we can prove the desired result. Under the conditions of the theorem, there exists $1+\sigma'_{i_0,j_0,k} g_{{i_0},{i_0}}^{\alpha}+\sigma''_{i_0,j_0,k} g_{{j_0},{j_0}}^{\beta}>0$.
Let $\! \mathbf{u}^k=[\mathbf{u}_0^k,\mathbf{u}_1^k,\ldots,\mathbf{u}_{N_x}^k]^T$,
where $ \mathbf{u}_i^k=[\epsilon_{i,0}^k,\epsilon_{i,1}^k,\ldots,\epsilon_{i,{N_y}}^k]^T$
and $||\mathbf{E}^k||_{\infty}=\max \limits _{{0\leq i \leq N_x},{0\leq j \leq N_y}}|\epsilon_{i,j}^k|$.
For  $k=0$, supposing  $|\epsilon_{i_0,j_0}^1|=||\mathbf{E}^1||_{\infty}=\max \limits _{{0\leq i \leq N_x},{0\leq j \leq N_y}}|\epsilon_{i,j}^1|$,
 we obtain
  \begin{align*}
&||\mathbf{E}^1||_{\infty}=|\epsilon_{i_0,j_0}^1| \\
& =|(1+\sigma'_{i_0,j_0,0} g_{{i_0},{i_0}}^{\alpha}+\sigma''_{i_0,j_0,0} g_{{j_0},{j_0}}^{\beta})\epsilon_{i_0,j_0}^{0}
                      +\sigma'_{i_0,j_0,0}\!\!\!\!\!\sum_{m=0,m\neq {i_0}}^{N_x}\!\!\!\!\!  g_{{i_0},m}^{\alpha}\epsilon_{m,j_0}^{0}
                       +\sigma''_{i_0,j_0,0}\!\!\!\!\!\sum_{m=0,m\neq {j_0}}^{N_y} \!\!\!\!\! g_{{j_0},m}^{\beta}\epsilon_{i_0,m}^{0}|\\
& \leq  (1+\sigma'_{i_0,j_0,0} g_{{i_0},{i_0}}^{\alpha}+\sigma''_{i_0,j_0,0} g_{{j_0},{j_0}}^{\beta})|\epsilon_{i_0,j_0}^{0}|
                      +\sigma'_{i_0,j_0,0}\!\!\!\!\!\sum_{m=0,m\neq {i_0}}^{N_x} \!\!\!\!\! g_{{i_0},m}^{\alpha}|\epsilon_{m,j_0}^{0}|
                      +\sigma''_{i_0,j_0,0}\!\!\!\!\!\sum_{m=0,m\neq {j_0}}^{N_y}\!\!\!\!\!  g_{{j_0},m}^{\beta}|\epsilon_{i_0,m}^{0}|\\
& \leq  (1+\sigma'_{i_0,j_0,0} g_{{i_0},{i_0}}^{\alpha}+\sigma''_{i_0,j_0,0}g_{{j_0},{j_0}}^{\beta})||\mathbf{E}^0||_{\infty}
                      +\sigma'_{i_0,j_0,0}\!\!\!\!\!\!\!\sum_{m=0,m\neq {i_0}}^{N_x} \!\!\!\!\!\! g_{{i_0},m}^{\alpha}||\mathbf{E}^0||_{\infty}
                      +\sigma''_{i_0,j_0,0}\!\!\!\!\!\!\!\sum_{m=0,m\neq {j_0}}^{N_y}\!\!\!\!\!\!  g_{{j_0},m}^{\beta}||\mathbf{E}^0||_{\infty}\\
&= ||\mathbf{E}^0||_{\infty}
                      +\sigma'_{i_0,j_0,0}\sum_{m=0}^{N_x}  g_{{i_0},m}^{\alpha}||\mathbf{E}^0||_{\infty}
                      +\sigma''_{i_0,j_0,0}\sum_{m=0}^{N_y}  g_{{j_0},m}^{\beta}||\mathbf{E}^0||_{\infty}
                      \leq ||\mathbf{E}^0||_{\infty}.
  \end{align*}
Assuming $||\mathbf{E}^{\widetilde k}||_{\infty} \leq ||\mathbf{E}^0||_{\infty}, {\widetilde k}=1,2,\ldots,k,$
and $|\epsilon_{i_0,j_0}^{k+1}|=||\mathbf{E}^{k+1}||_{\infty}=\max \limits _{{0\leq i \leq N_x},{0\leq j \leq N_y}}|\epsilon_{i,j}^{k+1}|$,
we get
  \begin{align*}
||\mathbf{E}^{k+1}||_{\infty}& =  |\epsilon_{i_0,j_0}^{k+1}|\\
 & =\Big|(1-l_1 +\sigma'_{i_0,j_0,k} g_{{i_0},{i_0}}^{\alpha}
                                                                       +\sigma''_{i_0,j_0,k} g_{{j_0},{j_0}}^{\beta})\epsilon_{i_0,j_0}^{k}
                            +\sum_{s=1}^{k-1}(l_s-l_{s+1})\epsilon_{i_0,j_0}^{k-s}+l_k\epsilon_{i_0,j_0}^{0}\\
                           & \quad+\sigma'_{i_0,j_0,k}\!\!\!\!\!\sum_{m=0,m\neq {i_0}}^{N_x} \!\!\!\!\! g_{{i_0},m}^{\alpha}\epsilon_{m,j_0}^{k}
                           +\sigma''_{i_0,j_0,k}\!\!\!\!\!\sum_{m=0,m\neq {j_0}}^{N_y}\!\!\!\!\!  g_{{j_0},m}^{\beta}\epsilon_{i_0,m}^{k}\Big|\\ 
                     &\leq (1-l_1 +\sigma'_{i_0,j_0,k} g_{{i_0},{i_0}}^{\alpha}
                            +\sigma''_{i_0,j_0,k} g_{{j_0},{j_0}}^{\beta})|\epsilon_{i_0,j_0}^{k}|
                            +\sum_{s=1}^{k-1}(l_s-l_{s+1})|\epsilon_{i_0,j_0}^{k-s}|+l_k|\epsilon_{i_0,j_0}^{0}|\\
                           & \quad+\sigma'_{i_0,j_0,k}\!\!\!\!\! \sum_{m=0,m\neq {i_0}}^{N_x} \!\!\!\!\! g_{{i_0},m}^{\alpha}|\epsilon_{m,j_0}^{k}|
                           +\sigma''_{i_0,j_0,k}\!\!\!\!\! \sum_{m=0,m\neq {j_0}}^{N_y}\!\!\!\!\!  g_{{j_0},m}^{\beta}|\epsilon_{i_0,m}^{k}|\\ 
                    &\leq (1-l_1)||\mathbf{E}^{k}||_{\infty}
                            +\sum_{s=1}^{k-1}(l_s-l_{s+1})||\mathbf{E}^{k-s}||_{\infty}+l_k||\mathbf{E}^{0}||_{\infty}\\
                    &      \quad  +\sigma'_{i_0,j_0,k}\sum_{m=0}^{N_x}  g_{{i_0},m}^{\alpha}||\mathbf{E}^{k}||_{\infty}
                            +\sigma''_{i_0,j_0,k}\sum_{m=0}^{N_y}  g_{{j_0},m}^{\beta}||\mathbf{E}^{k}||_{\infty}\\
                     &\leq (1-l_1)||\mathbf{E}^{0}||_{\infty}
                            +\sum_{s=1}^{k-1}(l_s-l_{s+1})||\mathbf{E}^{0}||_{\infty}+l_k||\mathbf{E}^{0}||_{\infty}\\
                     &= \Big(1-l_1+\sum_{s=1}^{k-1}(l_s-l_{s+1})+l_k\Big)||\mathbf{E}^{0}||_{\infty} =||\mathbf{E}^{0}||_{\infty}.
\end{align*}
When $\gamma=1$, it is easy to check that
\begin{align*}
&||\mathbf{E}^{k+1}||_{\infty} = |\epsilon_{i_0,j_0}^{k+1}|\\
                              &=\Big|(1+\sigma'_{i_0,j_0,k} g_{{i_0},{i_0}}^{\alpha}+\sigma''_{i_0,j_0,k} g_{{j_0},{j_0}}^{\beta})\epsilon_{i_0,j_0}^{k}
                                          +\sigma'_{i_0,j_0,k}\!\!\!\!\!\sum_{m=0,m\neq {i_0}}^{N_x} \!\!\!\!\!  g_{{i_0},m}^{\alpha}\epsilon_{m,j_0}^{k}
                                          +\sigma''_{i_0,j_0,k}\!\!\!\!\!\sum_{m=0,m\neq {j_0}}^{N_y}\!\!\!\!\!   g_{{j_0},m}^{\beta}\epsilon_{i_0,m}^{k}\Big|\\
                              &\leq(1+\sigma'_{i_0,j_0,k} g_{{i_0},{i_0}}^{\alpha}+\sigma''_{i_0,j_0,k} g_{{j_0},{j_0}}^{\beta})|\epsilon_{i_0,j_0}^{k}|
                                          +\sigma'_{i_0,j_0,k}\!\!\!\!\!\sum_{m=0,m\neq {i_0}}^{N_x} \!\!\!\!\!  g_{{i_0},m}^{\alpha}|\epsilon_{m,j_0}^{k}|
                                          +\sigma''_{i_0,j_0,k}\!\!\!\!\!\sum_{m=0,m\neq {j_0}}^{N_y} \!\!\!\!\!  g_{{j_0},m}^{\beta}|\epsilon_{i_0,m}^{k}|\\
                              &\leq (1+\sigma'_{i_0,j_0,k} g_{{i_0},{i_0}}^{\alpha}+\sigma''_{i_0,j_0,k} g_{{j_0},{j_0}}^{\beta})||\mathbf{E}^{k}||_{\infty}
                                         \! +\sigma'_{i_0,j_0,k}\!\!\!\!\!\!\sum_{m=0,m\neq {i_0}}^{N_x} \!\!\!\!\! \! g_{{i_0},m}^{\alpha}||\mathbf{E}^{k}||_{\infty}
                                          \!+\sigma''_{i_0,j_0,k}\!\!\!\!\!\!\sum_{m=0,m\neq {j_0}}^{N_y}\!\!\!\!\!\!   g_{{j_0},m}^{\beta}||\mathbf{E}^{k}||_{\infty}\\
                              &=||\mathbf{E}^{k}||_{\infty}+\sigma'_{i_0,j_0,k}\sum_{m=0}^{N_x}  g_{{i_0},m}^{\alpha}||\mathbf{E}^{k}||_{\infty}
                                        +\sigma''_{i_0,j_0,k}\sum_{m=0}^{N_y}  g_{{j_0},m}^{\beta}||\mathbf{E}^{k}||_{\infty}
                                          \leq ||\mathbf{E}^{k}||_{\infty}
                                          \leq ||\mathbf{E}^{0}||_{\infty}.
\end{align*}

\end{proof}

\section{Convergence analysis}
We use four subsections to prove that the global truncation error of the schemes (\ref{2.19}) and (\ref{2.24}) used to solve (\ref{2.16})  is $\mathcal{O} (\tau^{2-\gamma}+(\Delta x)^2)$, and the global truncation error of the schemes (\ref{2.22}) and (\ref{2.25}) used to solve (\ref{1.1}) is $\mathcal{O} (\tau^{2-\gamma}+(\Delta x)^2+(\Delta y)^2)$.


\subsection{Convergence for one-dimensional implicit difference scheme}

\noindent{\bf Theorem 4.1.}
\emph{ Let $u_i^k$ be the approximation solution of $u(x_i,t_k)$ computed by use of
the  implicit difference scheme (\ref{2.19}),  then there is a positive constant C such that
\begin{equation*}
  \begin{split}
|u(x_i,t_k)-u_i^k| \leq  C (\tau^{2-\gamma}+(\Delta x)^2), \quad i=0,1,\ldots,N_x;\,k=0,1,\ldots,N_t.
  \end{split}
  \end{equation*}.}
\begin{proof}
Let $u(x_i,t_k)$ be the exact solution of (\ref{2.16}) at the mesh point $(x_i,t_k)$.
Define $\varepsilon_i^k=u(x_i,t_k)-u_i^k$, and    $e^k=[\varepsilon_0^k,\varepsilon_1^k,\ldots, \varepsilon_{N_x}^k]$.
 Subtracting (\ref{2.17}) from (\ref{2.18}) and using $e^0=0$, we obtain
\begin{equation}\label{4.33}
  \begin{split}
(1-\omega_{i,1} g_{i,i}^{\alpha})\varepsilon_i^1-\omega_{i,1}\!\!\!\!\sum_{m=0,m\neq i}^{N_x} \!\!\!\!  g_{i,m}^{\alpha}\varepsilon_m^{1}&=R_i^{1}, ~~~k=0,\\
(1-\omega_{i,k+1} g_{i,i}^{\alpha})\varepsilon_i^{k+1}-\omega_{i,k+1}\!\!\!\!\sum_{m=0,m\neq i}^{N_x} \!\!\!\!  g_{i,m}^{\alpha}\varepsilon_m^{k+1}
  &=\sum_{s=0}^{k-1}(l_s-l_{s+1})\varepsilon_i^{k-s}+R_i^{k+1},    ~~k>0.
  \end{split}
  \end{equation}
When $\gamma=1$, (\ref{4.33}) can be written as
\begin{equation}\label{4.34}
  \begin{split}
(1-\omega_{i,1} g_{i,i}^{\alpha})\varepsilon_i^1-\omega_{i,1}\!\!\!\!\sum_{m=0,m\neq i}^{N_x} \!\!\!\!  g_{i,m}^{\alpha}\varepsilon_m^{1}&=R_i^{1}, ~~~k=0,\\
(1-\omega_{i,k+1} g_{i,i}^{\alpha})\varepsilon_i^{k+1}-\omega_{i,k+1}\!\!\!\!\sum_{m=0,m\neq i}^{N_x} \!\!\!\!  g_{i,m}^{\alpha}\varepsilon_m^{k+1}
  &=\varepsilon_i^{k}+R_i^{k+1},    ~~k>0.
  \end{split}
  \end{equation}
Denoting that  $||e^k||_{\infty}=\max \limits _{0\leq i \leq N_x}|\varepsilon_i^k|$
and $R_{\max}=\max \limits _{0\leq i \leq N_x,0\leq k \leq N_t}|R_i^k|$,  the desired result can be proved by using mathematical induction.

(1) Case $0 < \gamma < 1 $: 
For  $k=0$, supposing  $|\varepsilon_{i_0}^1|=||e^1||_{\infty}=\max \limits _{0\leq i \leq N_x}|\varepsilon_i^1|$,
according to Lemma 3.2, we get

  \begin{align*}
||e^1||_{\infty}&=|\varepsilon_{i_0}^1|  \leq |\varepsilon_{i_0}^1| -\omega_{i_0,1}\sum_{m=0}^{N_x}  g_{{i_0},m}^{\alpha}|\varepsilon_{i_0}^{1}|
                     =|\varepsilon_{i_0}^1|-\omega_{i_0,1} g_{{i_0},{i_0}}^{\alpha}|\varepsilon_{i_0}^{1}|
                                     -\omega_{i_0,1}\!\!\!\!\!\sum_{m=0,m\neq {i_0}}^{N_x} \!\!\!\! g_{{i_0},m}^{\alpha}|\varepsilon_{i_0}^{1}|\\
                 &\leq (1-\omega_{i_0,1} g_{{i_0},{i_0}}^{\alpha})|\varepsilon_{i_0}^{1}|
                                    -\omega_{i_0,1}\!\!\!\!\!\!\sum_{m=0,m\neq {i_0}}^{N_x} \!\!\!\!\! g_{{i_0},m}^{\alpha}|\varepsilon_m^{1}|
                    \leq \Big|(1-\omega_{i_0,1} g_{{i_0},{i_0}}^{\alpha})\varepsilon_{i_0}^{1}
                                    -\omega_{i_0,1}\!\!\!\!\!\sum_{m=0,m\neq {i_0}}^{N_x}\!\!\!\!  g_{{i_0},m}^{\alpha} \varepsilon_m^{1}\Big|\\
                 & =|R_{i_0}^1|  \leq R_{\max}=l_0^{-1}R_{\max}.
  \end{align*}
Supposing  $||e^{\widetilde k}||_{\infty} \leq l_{{\widetilde k}-1}^{-1}R_{\max}$, $ {\widetilde k}=1,2,\ldots,k,$
and $|\varepsilon_{i_0}^{k+1}|=||e^{k+1}||_{\infty}=\max \limits _{0\leq i \leq N_x}|\varepsilon_i^{k+1}|$,
according to the result of Lemma 3.1, $l_{{\widetilde k}}^{-1} \leq l_k^{-1}$, $ {\widetilde k}=0,1,\ldots,k$, therefore,
 $||e^{\widetilde k}||_{\infty} \leq l_{k}^{-1}R_{\max}$, $ {\widetilde k}=1,2,\ldots,k$.
Then we have 
  \begin{align*}
||e^{k+1}||_{\infty}&=|\varepsilon_{i_0}^{k+1}|  \leq |\varepsilon_{i_0}^{k+1}| -\omega_{i_0,k+1}\sum_{m=0}^{N_x}  g_{{i_0},m}^{\alpha}|\varepsilon_{i_0}^{k+1}|\\
                    &\quad  =|\varepsilon_{i_0}^{k+1}|- \omega_{i_0,k+1} g_{{i_0},{i_0}}^{\alpha}|\varepsilon_{i_0}^{k+1}|
                      -\omega_{i_0,k+1}\!\!\!\!\!\sum_{m=0,m\neq {i_0}}^{N_x} \!\!\!\!\! g_{{i_0},m}^{\alpha}|\varepsilon_{i_0}^{k+1}|\\
                    &\quad \leq (1-\omega_{i_0,k+1} g_{{i_0},{i_0}}^{\alpha})|\varepsilon_{i_0}^{k+1}| -\omega_{i_0,k+1}\!\!\!\!\!\sum_{m=0,m\neq {i_0}}^{N_x}  \!\!\!\!\! g_{{i_0},m}^{\alpha}|\varepsilon_m^{k+1}| \\
                     &\quad  \leq \Big|(1-\omega_{i_0,k+1} g_{{i_0},{i_0}}^{\alpha})\varepsilon_{i_0}^{k+1} -\omega_{i_0,k+1}\!\!\!\!\!\sum_{m=0,m\neq {i_0}}^{N_x}  \!\!\!\!\! g_{{i_0},m}^{\alpha} \varepsilon_m^{k+1}\Big|\\
                     &\quad=\Big|\sum_{s=0}^{k-1}(l_s-l_{s+1})\varepsilon_{i_0}^{k-s}+R_{i_0}^{k+1}\Big|\\
                    &\quad=\Big|(1-l_1)\varepsilon_{i_0}^{k} +\sum_{s=1}^{k-1}(l_s-l_{s+1})\varepsilon_{i_0}^{k-s}+R_{i_0}^{k+1}\Big|.
\end{align*}
From Lemma 3.1, we obtain

  \begin{align*}
||e^{k+1}||_{\infty} &\leq  (1-l_1) ||e^{k}||_{\infty}   +\sum_{s=1}^{k-1}(l_s-l_{s+1})||e^{k-s}||_{\infty}+|R_{i_0}^{k+1}|\\
                     & \leq l_k^{-1}\left (1-l_1+\sum_{s=1}^{k-1}(l_s-l_{s+1})+l_k\right)R_{\max} \\
                     & = l_k^{-1}R_{\max}  \leq Ck^{\gamma}R_{\max}
                     = C(k\tau)^{\gamma}(\tau^{2-\gamma}+(\Delta x)^2)\\
                     &\leq CT^{\gamma}(\tau^{2-\gamma}+(\Delta x)^2).
\end{align*}

(2) Case $\gamma =1$: 
Using similar idea leads to 
\begin{align*}
||e^{k+1}||_{\infty} &\leq  |(1-\omega_{i_0,k+1} g_{i_0,i_0}^{\alpha})\varepsilon_{i_0}^{k+1}
                  -\omega_{i_0,k+1}\!\!\!\!\!\sum_{m=0,m\neq i_0}^{N_x} \!\!\!\! g_{i_0,m}^{\alpha}\varepsilon_m^{k+1}|
  =|\varepsilon_{i_0}^{k}+R_{i_0}^{k+1}|\\
  &\leq ||e^{k}||_{\infty}+|R_{i_0}^{k+1}|\leq (k+1)R_{\max}\leq C(k+1)\tau(\tau+(\Delta x)^2)\\
  &\leq C'(\tau+(\Delta x)^2)=C'(\tau^{2-\gamma}+(\Delta x)^2).
\end{align*}
\end{proof}

\subsection{Convergence for one-dimensional explicit difference scheme}

\noindent{\bf Theorem 4.2.}
\emph{ Let $u_i^k$ be the approximation solution of $u(x_i,t_k)$ computed by use of
the  explicit difference scheme (\ref{2.25}). If
$0<\frac{\tau^{\gamma}}{(\Delta x)^{\alpha}} \leq -\frac{\Gamma(4-\alpha)(1-2^{-\gamma})}{4\kappa_{\alpha}C_{max}\Gamma(2-\gamma) (1-2^{1-\alpha})}$,
 where $C_{max}=\max\limits_{0\leq i \leq N_x,0\leq k \leq N_t} c(x_i,t_k),
$ then there is a positive constant C such that
\begin{equation*}
  \begin{split}
|u(x_i,t_k)-u_i^k| \leq  C (\tau^{2-\gamma}+(\Delta x)^2), \quad i=0,1,\ldots,N_x;\,k=0,1,\ldots,N_t.
  \end{split}
  \end{equation*}.}

\begin{proof}
Define $\varepsilon_i^k=u(x_i,t_k)-u_i^k$ and $R_{\max}=\max \limits _{0\leq i \leq N_x,0\leq k \leq N_t}|R_i^k|$.
Analogously, we have

\begin{equation}\label{4.35}
  \begin{split}
\epsilon_i^1&=R_i^1, \qquad k=0,\\
\varepsilon_i^{k+1}&=(1-l_1+\sigma_{i,k} g_{i,i}^{\alpha})\varepsilon_i^{k}+\sum_{s=1}^{k-1}(l_s-l_{s+1})\epsilon_i^{k-s}
+\sigma_{i,k}\!\!\!\!\!\sum_{m=0,m\neq i}^{N_x}\!\!\!\!\!  g_{i,m}^{\alpha}\varepsilon_m^{k}+R_i^{k+1},   \quad  k>0.
 \end{split}
\end{equation}
When $\gamma=1$, (\ref{4.35}) can be rewrote as 
\begin{equation}\label{4.36}
  \begin{split}
\epsilon_i^1&=R_i^1, \qquad k=0,\\
\varepsilon_i^{k+1}&=(1+\sigma_{i,k} g_{i,i}^{\alpha})\varepsilon_i^{k}
+\sigma_{i,k}\!\!\!\!\!\sum_{m=0,m\neq i}^{N_x}\!\!\!\!\!  g_{i,m}^{\alpha}\varepsilon_m^{k}+R_i^{k+1},   \qquad  k>0.
 \end{split}
\end{equation}
Denoting that $e^k=[\varepsilon_0^k,\varepsilon_1^k,\ldots, \varepsilon_{N_x}^k]$
and $||e^k||_{\infty}=\max \limits _{0\leq i \leq N_x}|\epsilon_i^k|$, we use mathematical induction to prove the desired result.

(1) Case $0 < \gamma < 1 $: 
For  $k=0$, supposing $|\varepsilon_{i_0}^1|=||e^1||_{\infty}=\max \limits _{0\leq i \leq N_x}|\varepsilon_i^1|$,
we have
  \begin{align*}
||e^1||_{\infty}=|\varepsilon_{i_0}^1|  = |R_{i_0}^1| \leq R_{\max} =l_0^{-1}R_{\max}.
  \end{align*}
Assuming  $||e^{\widetilde k}||_{\infty} \leq l_{{\widetilde k}-1}^{-1}R_{\max}$, $ {\widetilde k}=1,2,\ldots,k,$
and $|\varepsilon_{i_0}^{k+1}|=||e^{k+1}||_{\infty}=\max \limits _{0\leq i \leq N_x}|\varepsilon_i^{k+1}|$,
using  $l_{{\widetilde k}}^{-1} \leq l_k^{-1}$, $ {\widetilde k}=0,1,\ldots,k$, therefore,
 $||e^{\widetilde k}||_{\infty} \leq l_{k}^{-1}R_{\max}$, $ {\widetilde k}=1,2,\ldots,k$. Then we get  
  \begin{align*}
||e^{k+1}||_{\infty} &=  |\varepsilon_{i_0}^{k+1}|  =\Big|(1-\!l_1+\!\sigma_{i_0,k} g_{{i_0},{i_0}}^{\alpha})\varepsilon_{i_0}^{k}
                            +\sum_{s=1}^{k-1}(l_s-l_{s+1})\varepsilon_{i_0}^{k-s}
                            +\sigma_{i_0,k}\!\!\!\!\!\!\sum_{m=0,m\neq {i_0}}^{N_x}\!\!\!\!\!\!  g_{{i_0},m}^{\alpha}\varepsilon_m^{k}
                            +R_{i_0}^{k+1}\Big|\\
                     &\leq (1-l_1+\sigma_{i_0,k} g_{{i_0},{i_0}}^{\alpha})|\varepsilon_{i_0}^{k}|
                            +\sum_{s=1}^{k-1}(l_s-l_{s+1})|\varepsilon_{i_0}^{k-s}|
                            +\sigma_{i_0,k}\!\!\!\!\!\sum_{m=0,m\neq {i_0}}^{N_x} \!\!\!\!\! g_{{i_0},m}^{\alpha}|\varepsilon_m^{k}|
                            +R_{\max}\\
                     &\leq (1-l_1)||e^{k}||_{\infty}
                            +\sum_{s=1}^{k-1}(l_s-l_{s+1})||e^{k-s}||_{\infty}
                            +\sigma_{i_0,k}\sum_{m=0}^{N_x}  g_{{i_0},m}^{\alpha}||e^{k}||_{\infty}
                            +R_{\max}\\
                     &\leq l_k^{-1}\left(1-l_1
                            +\sum_{s=1}^{k-1}(l_s-l_{s+1})+l_k\right)R_{\max}\\
                     &=  l_k^{-1}R_{\max}\leq  C (\tau^{2-\gamma}+(\Delta x)^2).
\end{align*}

(2) Case $\gamma =1$:  Analogously, we have

\begin{align*}
||e^{k+1}||_{\infty} =|\varepsilon_{i_0}^{k+1}|
                                        &\leq (1+\sigma_{i_0,k} g_{{i_0},{i_0}}^{\alpha})|\varepsilon_{i_0}^{k}|
                                          +\sigma_{i_0,k}\!\!\!\!\!\sum_{m=0,m\neq {i_0}}^{N_x}  g_{{i_0},m}^{\alpha}|\varepsilon_m^{k}|+R_{\max}\\
                                         &\leq ||e^{k}||_{\infty}+\sigma_{i_0,k}\sum_{m=0}^{N_x}  g_{{i_0},m}^{\alpha} ||e^{k}||_{\infty}+R_{\max}\\
                                         &\leq ||e^{k}||_{\infty}+R_{\max} \leq C(\tau^{2-\gamma}+(\Delta x)^2).
\end{align*}
\end{proof}

\subsection{Convergence for two-dimensional implicit difference scheme}

\noindent{\bf Theorem 4.3.}
\emph{ Let $u_{i,j}^k$ be the approximation solution of $u(x_i,y_j,t_k)$ computed by use of
the  implicit difference scheme (\ref{2.22}),  then there is a positive constant C such that
\begin{equation*}
  \begin{split}
|u(x_i,y_j,t_k)-u_{i,j}^k| \leq  C (\tau^{2-\gamma}+(\Delta x)^2+(\Delta y)^2),
  \end{split}
  \end{equation*}
where $ i=0,1,\ldots,N_x;\,j=0,1,\ldots,N_y;\,k=0,1,\ldots,N_t.$}

\begin{proof}
Defining $\varepsilon_{i,j}^k=u(x_i,y_j,t_k)-u_{i,j}^k$,  $\mathbf{R}_{\max}=\max \limits _{0\leq i \leq N_x,0\leq j \leq N_y,0\leq k \leq N_t}|R_{i_j}^k|$, and subtracting (\ref{2.20}) from (\ref{2.21}), we obtain
\begin{equation}\label{4.37}
  \begin{split}
&(1-\omega'_{i,j,k+1} g_{i,i}^{\alpha}-\omega''_{i,j,k+1} g_{j,j}^{\beta})\varepsilon_{i,j}^{1}-\omega'_{i,j,k+1}\!\!\!\!\!\sum_{m=0,m\neq i}^{N_x}  \!\!\!\!\! g_{i,m}^{\alpha}\varepsilon_{m,j}^{1}
-\omega''_{i,j,k+1}\!\!\!\!\!\sum_{m=0,m\neq j}^{N_y} \!\!\!\!\! g_{j,m}^{\beta}\varepsilon_{i,m}^{1}\\
&\qquad\qquad=R_{i,j}^{1},  \quad k=0,\\
&(1-\omega'_{i,j,k+1}g_{i,i}^{\alpha}-\omega''_{i,j,k+1} g_{j,j}^{\beta})\varepsilon_{i,j}^{k+1}-\omega'_{i,j,k+1}\!\!\!\!\!\sum_{m=0,m\neq i}^{N_x}  \!\!\!\!\! g_{i,m}^{\alpha}\varepsilon_{m,j}^{k+1}
-\omega''_{i,j,k+1}\!\!\!\!\!\sum_{m=0,m\neq j}^{N_y}\!\!\!\!\!  g_{j,m}^{\beta}\varepsilon_{i,m}^{k+1}\\
 &\qquad\qquad=\sum_{s=0}^{k-1}(l_s-l_{s+1})\varepsilon_{i,j}^{k-s}+R_{i,j}^{k+1},    \quad k>0.
 \end{split}
\end{equation}
When $\gamma=1$, (\ref{4.37}) can be written as 
\begin{equation*}
\begin{split}
&(1-\omega'_{i,j,k+1} g_{i,i}^{\alpha}-\omega''_{i,j,k+1} g_{j,j}^{\beta})\varepsilon_{i,j}^{k+1}-\omega'_{i,j,k+1}\!\!\!\!\!\sum_{m=0,m\neq i}^{N_x}  \!\!\!\!\!g_{i,m}^{\alpha}\varepsilon_{m,j}^{k+1}
-\omega''_{i,j,k+1}\!\!\!\!\!\sum_{m=0,m\neq j}^{N_y} \!\!\!\!\! g_{j,m}^{\beta}\varepsilon_{i,m}^{k+1}\\
&\quad =\varepsilon_{i,j}^{k}+R_{i,j}^{k+1},       ~~ k \geq 0.
 \end{split}
\end{equation*}
Denoting $\mathbf{u}^k=[\mathbf{u}_0^k,\mathbf{u}_1^k,\ldots,\mathbf{u}_{N_x}^k]^T$,
 $\mathbf{u}_i^k=[\varepsilon_{i,0}^k,\varepsilon_{i,1}^k,\ldots,\varepsilon_{i,{N_y}}^k]^T$
and $||\mathbf{e}^k||_{\infty}=\max \limits _{{0\leq i \leq N_x},{0\leq j \leq N_y}}|\varepsilon_{i,j}^k|$, we prove the desired result by mathematical induction.

(1) Case $0 < \gamma < 1 $: 
For  $k=0$, supposing  $|\varepsilon_{i_0,j_0}^1|=||\mathbf{e}^1||_{\infty}=\max \limits _{{0\leq i \leq N_x},{0\leq j \leq N_y}}|\varepsilon_{i,j}^1|$,
 we have
  \begin{align*}
||\mathbf{e}^1||_{\infty}=& |\varepsilon_{i_0,j_0}^1|
                     \leq  |\varepsilon_{i_0,j_0}^1| -\omega'_{i_0,j_0,1}\sum_{m=0}^{N_x}  g_{{i_0},m}^{\alpha}|\varepsilon_{i_0,j_0}^{1}|
                            -\omega''_{i_0,j_0,1}\sum_{m=0}^{N_y}  g_{{j_0},m}^{\beta}|\varepsilon_{i_0,j_0}^{1}|\\
                      \leq & \Big |(1-\omega'_{i_0,j_0,1} g_{{i_0},{i_0}}^{\alpha}-\omega''_{i_0,j_0,1} g_{{j_0},{j_0}}^{\beta})\varepsilon_{i_0,j_0}^{1}\\
                           & -\omega'_{i_0,j_0,1}\!\!\!\!\!\sum_{m=0,m\neq {i_0}}^{N_x}  g_{{i_0},m}^{\alpha}\varepsilon_{m,{j_0}}^{1}
                            -\omega''_{i_0,j_0,1}\!\!\!\!\!\sum_{m=0,m\neq {j_0}}^{N_y}  g_{{j_0},m}^{\beta}\varepsilon_{{i_0,m}}^{1}\Big|\\
                          =&|R_{i_0,j_0}^{1}| \leq \mathbf{R}_{\max} =l_0^{-1}\mathbf{R}_{\max}.
  \end{align*}
Assuming $||e^{\widetilde k}||_{\infty} \leq l_{{\widetilde k}-1}^{-1}\mathbf{R}_{\max}$, $ {\widetilde k}=1,2,\ldots,k,$
and $|\varepsilon _{i_0,j_0}^{k+1}|=||\mathbf{e}^{k+1}||_{\infty}=\max \limits _{{0\leq i \leq N_x},{0\leq j \leq N_y}}|\epsilon_{i,j}^{k+1}|$,
using  $l_{{\widetilde k}}^{-1} \leq l_k^{-1}$, $ {\widetilde k}=0,1,\ldots,k$, therefore,
 $||\mathbf{e}^{\widetilde k}||_{\infty} \leq l_{k}^{-1}\mathbf{R}_{\max}$, $ {\widetilde k}=1,2,\ldots,k$.
Then we have 
  \begin{align*}
||\mathbf{e}^{k+1}||_{\infty}&= |\varepsilon_{i_0,j_0}^{k+1}|  \leq |\varepsilon_{i_0,j_0}^{k+1}| -\omega'_{i_0,j_0,k+1}\sum_{m=0}^{N_x}  g_{{i_0},m}^{\alpha}|\epsilon_{i_0,j_0}^{k+1}|
                            -\omega''_{i_0,j_0,k+1}\sum_{m=0}^{N_y}  g_{{j_0},m}^{\beta}|\varepsilon_{i_0,j_0}^{k+1}|\\
                          &~\quad  -\omega''_{i_0,j_0,k+1} g_{{j_0},{j_0}}^{\beta}|\epsilon_{i_0,j_0}^{k+1}|
                            -\omega''_{i_0,j_0,k+1}\!\!\!\!\!\sum_{m=0,m\neq {j_0}}^{N_y}\!\!\!\!\! g_{{j_0},m}^{\beta}|\varepsilon_{i_0,j_0}^{k+1}|\\
                     &~ \leq  (1-\omega'_{i_0,j_0,k+1} g_{{i_0},{i_0}}^{\alpha}-\omega''_{i_0,j_0,k+1} g_{{j_0},{j_0}}^{\beta})|\varepsilon_{i_0,j_0}^{k+1}|\\
                     &~\quad      -\omega'_{i_0,j_0,k+1}\!\!\!\!\!\sum_{m=0,m\neq {i_0}}^{N_x}  g_{{i_0},m}^{\alpha}|\varepsilon_{m,{j_0}}^{k+1}|
                            -\omega''_{i_0,j_0,k+1}\!\!\!\!\!\sum_{m=0,m\neq {j_0}}^{N_y}  g_{{j_0},m}^{\beta}|\varepsilon_{{i_0,m}}^{k+1}|\\
                     &~ \leq \Big |(1-\omega'_{i_0,j_0,k+1} g_{{i_0},{i_0}}^{\alpha}-\omega''_{i_0,j_0,k+1} g_{{j_0},{j_0}}^{\beta})\varepsilon_{i_0,j_0}^{k+1}\\
                     &~\quad -\omega'_{i_0,j_0,k+1}\!\!\!\!\!\sum_{m=0,m\neq {i_0}}^{N_x}  g_{{i_0},m}^{\alpha}\varepsilon_{m,{j_0}}^{k+1}
                            -\omega''_{i_0,j_0,k+1}\!\!\!\!\!\sum_{m=0,m\neq {j_0}}^{N_y}  g_{{j_0},m}^{\beta}\varepsilon_{{i_0,m}}^{k+1}\Big|\\
                    &~=\big|\sum_{s=0}^{k-1}(l_s-l_{s+1})\varepsilon_{i_0,j_0}^{k-s}+R_{i_0,j_0}^{k+1}\big|\\
                    &\leq \sum_{s=0}^{k-1}(l_s-l_{s+1})||\mathbf{e}^{k-s}||_{\infty}
                            +\mathbf{R}_{\max}\\
                     &\leq l_k^{-1}\left(1-l_1
                            +\sum_{s=1}^{k-1}(l_s-l_{s+1})+l_k\right)\mathbf{R}_{\max}\\
                     &=  l_k^{-1}\mathbf{R}_{\max}\leq  C (\tau^{2-\gamma}+(\Delta x)^2+(\Delta y)^2).
\end{align*}

(2) Case $\gamma =1$: 
Using similar idea leads to 
\begin{align*}
||\mathbf{e}^{k+1}||_{\infty} &\leq \Big|(1-\omega'_{i_0,j_0,k+1} g_{i_0,i_0}^{\alpha}-\omega''_{i_0,j_0,k+1} g_{j_0,j_0}^{\beta})\varepsilon_{i_0,j_0}^{k+1}\\
     &~\quad  -\omega'_{i_0,j_0,k+1}\!\!\!\!\!\sum_{m=0,m\neq i}^{N_x}  g_{i_0,m}^{\alpha}\varepsilon_{m,j_0}^{k+1}
     -\omega''_{i_0,j_0,k+1}\!\!\!\!\!\sum_{m=0,m\neq j}^{N_y}  g_{j_0,m}^{\beta}\varepsilon_{i_0,m}^{k+1}\Big|\\
 &=|\varepsilon_{i_0,j_0}^{k}+R_{i_0,j_0}^{k+1}|
   \leq ||\mathbf{e}^{k}||_{\infty}+|R_{i_0,j_0}^{k+1}|\leq (k+1)\mathbf{R}_{\max}\\
  &\leq C(\tau^{2-\gamma}+(\Delta x)^2+(\Delta y)^2).
 \end{align*}
\end{proof}

\subsection{Convergence for two-dimensional explicit difference scheme}

\noindent{\bf Theorem 4.4.}
\emph{ Let $u_{i,j}^k$ be the approximation solution of $u(x_i,y_j,t_k)$ computed by use of
the  implicit difference scheme (\ref{2.25}). If
$
0<\frac{\tau^{\gamma}}{(\Delta x)^{\alpha}} +\frac{\tau^{\gamma}}{(\Delta y)^{\beta}}\leq \frac{1-2^{-\gamma}}{\Gamma(2-\gamma)\mathbf{C}_{max} }$,
 where $$\mathbf{C}_{max}=\max_{0\leq i \leq N_x,0\leq j \leq N_y,0\leq k \leq N_t}
\left \{-\frac{\kappa_{\alpha}(4-2^{3-\alpha})c_{i,j}^k}{\Gamma (4-\alpha)},-\frac{\kappa_{\beta}(4-2^{3-\beta})d_{i,j}^k}{\Gamma (4-\beta)}\right\},
$$  then there is a positive constant C such that
\begin{equation*}
  \begin{split}
|u(x_i,y_j,t_k)-u_{i,j}^k| \leq  C (\tau^{2-\gamma}+(\Delta x)^2+(\Delta y)^2),
  \end{split}
  \end{equation*}
where $ i=0,1,\ldots,N_x;\,j=0,1,\ldots,N_y;\,k=0,1,\ldots,N_t.$}

\begin{proof}
Define $\varepsilon_{i,j}^k=u(x_i,y_j,t_k)-u_{i,j}^k$ and  $\mathbf{R}_{\max}=\max \limits _{0\leq i \leq N_x,0\leq j \leq N_y,0\leq k \leq N_t}|R_{i,j}^k|$.
Similarly,   we obtain the following  form
\begin{equation}\label{4.41}
  \begin{split}
\varepsilon_{i,j}^{1}&=R_{i,j}^1, \quad k=0,\\
\varepsilon_{i,j}^{k+1}&=(1-l_1+\sigma'_{i,j,k} g_{i,i}^{\alpha}+\sigma''_{i,j,k} g_{j,j}^{\beta})\varepsilon_{i,j}^{k}+\sum_{s=1}^{k-1}(l_s-l_{s+1})\varepsilon_{i,j}^{k-s}\\
&\quad+\sigma'_{i,j,k}\!\!\!\!\!\sum_{m=0,m\neq i}^{N_x}  g_{i,m}^{\alpha}\varepsilon_{m,j}^{k}
       +\sigma''_{i,j,k}\!\!\!\!\!\sum_{m=0,m\neq j}^{N_y}  g_{j,m}^{\beta}\varepsilon_{i,m}^{k}+R_{i,j}^{k+1},    \quad k>0.
  \end{split}
\end{equation}
When  $\gamma=1$, (\ref{4.41}) becomes  
\begin{align*}
\varepsilon_{i,j}^{1}&=R_{i,j}^1, \quad k=0,\\
\epsilon_{i,j}^{k+1}&=(1+\sigma'_{i,j,k} g_{i,i}^{\alpha}+\sigma''_{i,j,k} g_{j,j}^{\beta})\varepsilon_{i,j}^{k}\\
&\quad +\sigma'_{i,j,k}\!\!\!\!\!\sum_{m=0,m\neq i}^{N_x}  g_{i,m}^{\alpha}\varepsilon_{m,j}^{k}
+\sigma''_{i,j,k}\!\!\!\!\!\sum_{m=0,m\neq j}^{N_y}  g_{j,m}^{\beta}\varepsilon_{i,m}^{k}+R_{i,j}^{k+1}, \quad k \geq 0.
\end{align*}
Using  mathematical induction, we can obtain the desired result. 
Let $\! \mathbf{u}^k=[\mathbf{u}_0^k,\mathbf{u}_1^k,\ldots,\mathbf{u}_{N_x}^k]^T$,
where $ \mathbf{u}_i^k=[\varepsilon_{i,0}^k,\varepsilon_{i,1}^k,\ldots,\varepsilon_{i,{N_y}}^k]^T$
and $||\mathbf{e}^k||_{\infty}=\max \limits _{{0\leq i \leq N_x},{0\leq j \leq N_y}}|\varepsilon_{i,j}^k|$.

(1) Case $0 < \gamma < 1 $: 
For  $k=0$, supposing $|\varepsilon_{i_0,j_0}^1|=||\mathbf{e}^1||_{\infty}=\max \limits _{{0\leq i \leq N_x},{0\leq j \leq N_y}}|\varepsilon_{i,j}^1|$,
 we obtain
  \begin{align*}
&||\mathbf{e}^1||_{\infty}=|\varepsilon_{i_0,j_0}^1|=|R_{i_0,j_0}^1| \leq \mathbf{R}_{\max} =l_0^{-1}\mathbf{R}_{\max}.
  \end{align*}
Assuming $||e^{\widetilde k}||_{\infty} \leq l_{{\widetilde k}-1}^{-1}\mathbf{R}_{\max}$, $ {\widetilde k}=1,2,\ldots,k,$
and $|\varepsilon _{i_0,j_0}^{k+1}|=||\mathbf{e}^{k+1}||_{\infty}=\max \limits _{{0\leq i \leq N_x},{0\leq j \leq N_y}}|\epsilon_{i,j}^{k+1}|$,
using  $l_{{\widetilde k}}^{-1} \leq l_k^{-1}$, $ {\widetilde k}=0,1,\ldots,k$, therefore,
 $||\mathbf{e}^{\widetilde k}||_{\infty} \leq l_{k}^{-1}\mathbf{R}_{\max}$, $ {\widetilde k}=1,2,\ldots,k$.
Then we have 
  \begin{align*}
||\mathbf{e}^{k+1}||_{\infty}& =  |\varepsilon_{i_0,j_0}^{k+1}|\\
 & =\Big|(1-l_1 +\sigma'_{i_0,j_0,k} g_{{i_0},{i_0}}^{\alpha}
                                                                       +\sigma''_{i_0,j_0,k} g_{{j_0},{j_0}}^{\beta})\varepsilon_{i_0,j_0}^{k}
                            +\sum_{s=1}^{k-1}(l_s-l_{s+1})\varepsilon_{i_0,j_0}^{k-s}\\
                           & \quad+\sigma'_{i_0,j_0,k}\!\!\!\!\!\sum_{m=0,m\neq {i_0}}^{N_x} \!\!\!\!\! g_{{i_0},m}^{\alpha}\varepsilon_{m,j_0}^{k}
                           +\sigma''_{i_0,j_0,k}\!\!\!\!\!\sum_{m=0,m\neq {j_0}}^{N_y}\!\!\!\!\!  g_{{j_0},m}^{\beta}\varepsilon_{i_0,m}^{k}+R_{i,j}^{k+1}\Big|\\ 
                     &\leq (1-l_1 +\sigma'_{i_0,j_0,k} g_{{i_0},{i_0}}^{\alpha}
                            +\sigma''_{i_0,j_0,k} g_{{j_0},{j_0}}^{\beta})|\varepsilon_{i_0,j_0}^{k}|
                            +\sum_{s=1}^{k-1}(l_s-l_{s+1})|\varepsilon_{i_0,j_0}^{k-s}|\\
                           & \quad+\sigma'_{i_0,j_0,k}\!\!\!\!\! \sum_{m=0,m\neq {i_0}}^{N_x} \!\!\!\!\! g_{{i_0},m}^{\alpha}|\varepsilon_{m,j_0}^{k}|
                           +\sigma''_{i_0,j_0,k}\!\!\!\!\! \sum_{m=0,m\neq {j_0}}^{N_y}\!\!\!\!\!  g_{{j_0},m}^{\beta}|\varepsilon_{i_0,m}^{k}|+\mathbf{R}_{\max}\\ 
                    &\leq (1-l_1)||\mathbf{e}^{k}||_{\infty}
                            +\sum_{s=1}^{k-1}(l_s-l_{s+1})||\mathbf{e}^{k-s}||_{\infty}+\\
                    &      \quad  +\sigma'_{i_0,j_0,k}\sum_{m=0}^{N_x}  g_{{i_0},m}^{\alpha}||\mathbf{e}^{k}||_{\infty}
                            +\sigma''_{i_0,j_0,k}\sum_{m=0}^{N_y}  g_{{j_0},m}^{\beta}||\mathbf{e}^{k}||_{\infty}+\mathbf{R}_{\max}\\
                    &\leq l_k^{-1}\left(1-l_1
                            +\sum_{s=1}^{k-1}(l_s-l_{s+1})+l_k\right)\mathbf{R}_{\max}\\
                     &=  l_k^{-1}\mathbf{R}_{\max}\leq  C (\tau^{2-\gamma}+(\Delta x)^2+(\Delta y)^2).
\end{align*}

(2) Case $\gamma =1$: It is easy to check that
\begin{align*}
||\mathbf{e}^{k+1}||_{\infty}& = |\varepsilon_{i_0,j_0}^{k+1}|\\
                              &\leq(1+\sigma'_{i_0,j_0,k} g_{{i_0},{i_0}}^{\alpha}+\sigma''_{i_0,j_0,k} g_{{j_0},{j_0}}^{\beta})|\varepsilon_{i_0,j_0}^{k}|\\
                               & \quad   +\sigma'_{i_0,j_0,k}\!\!\!\!\!\sum_{m=0,m\neq {i_0}}^{N_x} \!\!\!\!\!  g_{{i_0},m}^{\alpha}|\varepsilon_{m,j_0}^{k}|
                                          +\sigma''_{i_0,j_0,k}\!\!\!\!\!\sum_{m=0,m\neq {j_0}}^{N_y} \!\!\!\!\!  g_{{j_0},m}^{\beta}|\varepsilon_{i_0,m}^{k}|
                                          +\mathbf{R}_{\max}\\
                              &\leq (1+\sigma'_{i_0,j_0,k} g_{{i_0},{i_0}}^{\alpha}+\sigma''_{i_0,j_0,k} g_{{j_0},{j_0}}^{\beta})||\mathbf{e}^{k}||_{\infty}\\
                              &\quad +\sigma'_{i_0,j_0,k}\!\!\!\!\!\!\sum_{m=0,m\neq {i_0}}^{N_x} \!\!\!\!\! \! g_{{i_0},m}^{\alpha}||\mathbf{e}^{k}||_{\infty}
                                          \!+\sigma''_{i_0,j_0,k}\!\!\!\!\!\!\sum_{m=0,m\neq {j_0}}^{N_y}\!\!\!\!\!\!   g_{{j_0},m}^{\beta}||\mathbf{e}^{k}||_{\infty}+\mathbf{R}_{\max}\\
                                          &\leq ||\mathbf{e}^{k}||_{\infty}+\mathbf{R}_{\max}\leq (k+1)\mathbf{R}_{\max}\\
  &\leq C(\tau^{2-\gamma}+(\Delta x)^2+(\Delta y)^2).
\end{align*}

\end{proof}

\noindent{\bf Remark 4.4.} \label{Remark:12}
Let $\nu \in (0,1)$, there exist $_{x_L}\!D_x^{\nu} u=D{_{x_L}\!D_x^{-(1-\nu)}}u=D^2{_{x_L}\!D_x^{-(2-\nu)}}u$, and 
 $_x\!D_{x_R}^{\nu} u=-D{_x}\!D_{x_R}^{-(1-\nu)}u=D^2{_x\!D_{x_R}^{-(2-\nu)}}u$ \cite{Miller:93}. So after carefully dealing with the boundary conditions, all the above presented numerical schemes still work well when the order of space fractional derivatives $\alpha \in (0,1)$ and/or $\beta\in (0,1)$. And the convergence rates remain, and all the theoretical analyses are valid.  Some of the numerical results are also given in Table \ref{table:b}.


\section{Numerical results}
In this section, we numerically verify the above theoretical results including convergence
rates and numerical stability.  And the $ l_\infty$ norm is used to measure the numerical errors. 

%

\subsection{Numerical results for the implicit scheme for the one-dimensional time-space Caputo-Riesz fractional diffusion equation}
Consider  the one-dimensional time-space Caputo-Riesz fractional   diffusion equation (\ref{2.16}),
on a finite domain  $0\leq x \leq 1 $,  $0<t \leq 1/2$,  with the coefficient $c(x,t)=x^{\alpha}t^{1-\gamma}$,
 the forcing function
\begin{equation*}
\begin{split}
 f(x,t)=&\frac{1}{2}\Gamma(3+\gamma)t^2x^2(x-1)^2 \\
 &  +\frac{t^3x^{\alpha}}{cos(\alpha \pi/2)}
   \left[ \frac{x^{2-\alpha}+(1-x)^{2-\alpha}}{\Gamma(3-\alpha)}   -6\frac{x^{3-\alpha}+(1-x)^{3-\alpha}}{\Gamma(4-\alpha)}
   +12\frac{x^{4-\alpha}+(1-x)^{4-\alpha}}{\Gamma(5-\alpha)} \right], 
  \end{split}
\end{equation*}
the initial condition $u(x,0)=0 $, and the boundary
conditions $u(0,t)=u(1,t)=0$. This fractional diffusion equation has the exact
value $u(x,t)=t^{2+\gamma}x^2(1-x)^2$, which may be confirmed by applying
the fractional differential equations
\begin{equation*}
\begin{split}
&_{x_L}D_{x}^{\nu}(x-x_L)^{\upsilon}=\frac{\Gamma({\upsilon}+1)}{\Gamma({\upsilon}+1-\nu)}(x-x_L)^{{\upsilon}-\nu},\\
&_xD_{x_R}^{\nu}(x_R-x)^{\upsilon}=\frac{\Gamma({\upsilon}+1)}{\Gamma({\upsilon}+1-\nu)}(x_R-x)^{{\upsilon}-\nu}.
\end{split}
\end{equation*}

\begin{table}[h]\fontsize{9.5pt}{12pt}\selectfont
  \begin{center}
  \caption{The maximum errors and convergence rates for the implicit scheme (\ref{2.19}) of
the  one-dimensional time-space Caputo-Riesz fractional diffusion equation (\ref{2.16}) with variable coefficient $c(x,t)=x^{\alpha}t^{1-\gamma}$ at $t=1/2$, and the time and space stepsizes are equal, i.e, $\tau=\Delta x$. }\vspace{5pt}
    \begin{tabular*}{\linewidth}{@{\extracolsep{\fill}}*{7}{c|c|c|c|c|c|c}}                                  \hline  
$\tau,\Delta x$ & $\alpha=1.2,\gamma=0.9$ & Rate   & $\alpha=1.2,\gamma=0.5$ & Rate   & $\alpha=1.2,\gamma=0.1$ &   Rate     \\\hline
              ~~1/40&             3.1438e-004&        &             6.3187e-005&        &             2.7395e-005&          \\\hline
              ~~1/80&             1.4713e-004&  1.0954&             2.2183e-005&  1.5102&             7.6703e-006&  1.8366   \\\hline
             ~~1/160&             6.8748e-005&  1.0977&             7.7888e-006&  1.5100&             2.0317e-006&  1.9166   \\\hline
             ~~1/320&             3.2097e-005&  1.0989&             2.7378e-006&  1.5084&             5.2912e-007&  1.9410   \\\hline
$ $ & $\alpha=1.9,\gamma=0.9$ & Rate   & $\alpha=1.9,\gamma=0.5$ & Rate   & $\alpha=1.9,\gamma=0.1$ &   Rate     \\\hline
              ~~1/40&             2.7655e-004&        &             5.6774e-005&        &             2.1114e-005&          \\\hline
              ~~1/80&             1.2919e-004&  1.0981&             1.9669e-005&  1.5293&             5.4717e-006&  1.9482   \\\hline
             ~~1/160&             6.0294e-005&  1.0994&             6.8351e-006&  1.5249&             1.4145e-006&  1.9517   \\\hline
             ~~1/320&             2.8133e-005&  1.0997&             2.3841e-006&  1.5195&             3.6518e-007&  1.9536   \\\hline
    \end{tabular*}\label{tab:a}
  \end{center}
\end{table}
Table \ref{tab:a} shows the maximum errors, at time $t=1/2$ with 
$\tau=\Delta x $, between the exact analytical values and the
numerical values obtained by applying the implicit scheme
(\ref{2.19}). Since the scheme has the global truncation error $\mathcal{O} (\tau^{2-\gamma}+(\Delta x)^2)$, the convergent rate should be $2-\gamma$ being confirmed by the numerical results.

\begin{table}[h]\fontsize{9.5pt}{12pt}\selectfont
  \begin{center}
  \caption{The maximum errors and convergent rates for the implicit scheme (\ref{2.19}) of
the one-dimensional time-space Caputo-Riesz fractional diffusion equation (\ref{2.16}) at $t=1/2$ with variable coefficient $c(x,t)=x^{\alpha}t^{1-\gamma}$, and $\tau=(\Delta x)^{\frac{2}{2-\gamma}}$, where $\gamma=0.9$.}\vspace{5pt}
    \begin{tabular*}{\linewidth}{@{\extracolsep{\fill}}*{9}{c|c|c|c|c|c|c|c|c}}                                   \hline  
          $\tau,\Delta x$&  $\alpha=1.9$&   Rate &   $\alpha=1.5$&   Rate&   $\alpha=1.2$& Rate&   $\alpha=0.3$&  Rate     \\\hline
              ~ 1/10&   2.4699e-004&        &   2.5801e-004&        &   2.5510e-004&          &2.4135e-004&            \\\hline
              ~ 1/20&   6.2966e-005&  1.9718&   6.5569e-005&  1.9763&   6.4560e-005&  1.9823  &6.1043e-005&  1.9832     \\\hline
              ~ 1/40&   1.5697e-005&  2.0041&   1.6226e-005&  2.0148&   1.5934e-005&  2.0185  &1.5085e-005&  2.0167     \\\hline
              ~ 1/80&   3.9368e-006&  1.9954&   4.0475e-006&  2.0032&   4.0433e-006&  1.9785  &3.7583e-006&  2.0050     \\\hline
    \end{tabular*}\label{table:b} \vspace{20pt}
  \end{center}
\end{table}
 Table \ref{table:b} shows the maximum errors at time $t=1/2$, and the time and space stepsizes are taken as 
$\tau=(\Delta x)^{\frac{2}{2-\gamma}}$. The numerical results confirm the 2nd order convergence in space directions. In particular, the numerical results when $\alpha=0.3\in(0,1)$ are also presented, which confirm the statement of Remark \ref{Remark:12}.


\subsection{Implicit schemes results for  two-dimensional  time-space Caputo-Riesz fractional    diffusion equation}

We further examine the two-dimensional time-space Caputo-Riesz fractional diffusion equation
 (\ref{1.1}),
on a finite domain  $0\leq x \leq 1 $,  $0<t \leq 1/2$,  with the variable coefficients
$c(x,y,t)=2x^{\alpha}y^{\beta}t^{1-\gamma}$, $d(x,y,t)=2x^{\beta}y^{\alpha}t^{1-\gamma}$,
 the forcing function
\begin{align*}
 &f(x,y,t)\\
 &=\frac{1}{2}\Gamma(3+\gamma)t^2x^2(x-1)^2y^2(y-1)^2 \\
 & \quad +\frac{t^3x^{\alpha}y^{2+\beta}(y-1)^2}{cos(\alpha \pi/2)}
   \left[ \frac{x^{2-\alpha}+(1-x)^{2-\alpha}}{\Gamma(3-\alpha)}   -6\frac{x^{3-\alpha}+(1-x)^{3-\alpha}}{\Gamma(4-\alpha)}
   +12\frac{x^{4-\alpha}+(1-x)^{4-\alpha}}{\Gamma(5-\alpha)} \right]    \\
 & \quad +\frac{t^3x^{2+\beta}(x-1)^2y^{\alpha}}{cos(\beta \pi/2)}
   \left[ \frac{y^{2-\beta}+(1-y)^{2-\beta}}{\Gamma(3-\beta)}   -6\frac{y^{3-\beta}+(1-y)^{3-\beta}}{\Gamma(4-\beta)}
   +12\frac{y^{4-\beta}+(1-y)^{4-\beta}}{\Gamma(5-\beta)} \right],
  \end{align*}
the initial condition $u_0(x,y)=0 $, and the boundary
conditions $B(x,y,t)=0$. It has the exact solution 
\begin{equation*}
  u(x,y,t)=t^{2+\gamma}x^2(1-x)^2y^2(1-y)^2.
\end{equation*}

\begin{table}[h]\fontsize{9.5pt}{12pt}\selectfont
  \begin{center}
  \caption{The maximum errors and convergence rates for the implicit scheme (\ref{2.22}) of
the  two-dimensional time-space Caputo-Riesz fractional diffusion equation (\ref{1.1}) at $t=1/2$ and $\tau=\Delta x=\Delta y$,
 with variable coefficients $c(x,y,t)=2x^{\alpha}y^{\beta}t^{1-\gamma}$
 and $d(x,y,t)=2x^{\beta}y^{\alpha}t^{1-\gamma}$.}\vspace{5pt}

    \begin{tabular*}{\linewidth}{@{\extracolsep{\fill}}*{7}{c|c|c|c|c|c|c}}                                    \hline 
 $\tau,\Delta x,\Delta y$ & $\gamma=0.9$ &    & $\gamma=0.5$ &    & $\gamma=0.1$ &    \\\hline
 $  $ & $\alpha=1.2,\beta=1.3$ & Rate   & $\alpha=1.2,\beta=1.3$ & Rate   & $\alpha=1.2,\beta=1.3$ &  Rate  \\\hline
 ~~1/10&  7.7867E-005&        &  3.1936E-005&        &  2.2077E-005&          \\\hline
 ~~1/20&  3.6381E-005&  1.0978&  1.0357E-005&  1.6246&  5.6507E-006 &  1.9660   \\\hline
 ~~1/30&  2.3084E-005&  1.1220&  5.4951E-006&  1.5632&  2.5497E-006& 1.9627   \\\hline
 ~~1/40&  1.6839E-005&  1.0965&  3.4925E-006&  1.5755&  1.4509E-006 &  1.9598   \\\hline
 $$ & $\alpha=1.8,\beta=1.7$ & Rate   & $\alpha=1.8,\beta=1.7$ & Rate   & $\alpha=1.8,\beta=1.7$ &  Rate  \\\hline
 ~~1/10&  7.8209E-005&        &  3.1251E-005&        &  2.0743E-005&          \\\hline
 ~~1/20&  3.6912E-005&  1.0832&  1.0369E-005&  1.5917&  5.5573E-006 &  1.9002   \\\hline
 ~~1/30&  2.3535E-005&  1.1099&  5.5386E-006&  1.5465&  2.5115E-006 &  1.9588   \\\hline
 ~~1/40&  1.7122E-005&  1.1060&  3.5433E-006&  1.5527&  1.4310E-006 &  1.9552   \\\hline
    \end{tabular*}\label{tab:c}
  \end{center}
\end{table}
Table \ref{tab:c} shows the maximum errors of the implicit scheme
(\ref{2.22}) at time $t=1/2$, and
$\tau=\Delta x=\Delta y $. Since the implicit scheme
(\ref{2.22}) has the global truncation error $\mathcal{O} (\tau^{2-\gamma}+(\Delta x)^2+(\Delta y)^2)$, the convergence rate should be $2-\gamma$ being confirmed by the numerical results.

\section{Conclusions}
This paper discusses the finite difference schemes for the time-space Caputo-Riesz fractional PDEs with variable coefficients, and the order of time fractional derivative belongs to $(0,1)$ and the order of space fractional derivatives locate in $(1,2)$ . The obtained schemes have $(2-\gamma)$-th order convergence rate in time and 2nd order convergent rate in space. The detailed numerical stability analysis and error estimates are presented, and the extensive numerical experiments are performed, which confirm the theoretical results. In particular, the numerical schemes still work well for the time-space Caputo-Riesz fractional PDEs with the order of its space derivatives belongs to $(0,1)$, and all the theoretical analyses are still valid.


\section*{Acknowledgments} This work was supported by the Program for
New Century Excellent Talents in University under Grant No.
NCET-09-0438, the National Natural Science Foundation of China under
Grant No. 10801067, and the Fundamental Research Funds for the
Central Universities under Grant No. lzujbky-2010-63 and No. lzujbky-2012-k26.

\end{document}